\documentclass[10pt]{article}

\usepackage{latexsym,amssymb,amsmath,theorem,tan-gle,makeidx}
\usepackage{diagrams,hyperref}
  {\bigskip\noindent\begin{Sbox}\begin{minipage}{0.95\textwidth} }%
  {\end{minipage}\end{Sbox}\fbox{\TheSbox}\bigskip}%

\newlength{\dinwidth}
\newlength{\dinmargin}
\def\endexem{\hfill{$\Box$}\medskip}

\theoremstyle{change}
\theoremheaderfont{\scshape}
\newtheorem{thm}{Theorem}[section]
\newtheorem{prop}[thm]{Proposition}
\newtheorem{lem}[thm]{Lemma}
\newtheorem{cor}[thm]{Corollary}
\newtheorem{conj}[thm]{Conjecture}
\newtheorem{defin}[thm]{Definition}

\newtheorem{defprop}[thm]{Definition/Proposition}
\theorembodyfont{\rmfamily}
\newtheorem{leer}[thm]{}
\newtheorem{example}[thm]{Example}
\newtheorem{rema}[thm]{Remark}

\newcommand{\be}{\begin{equation}}
\newcommand{\ee}{\end{equation}}
\newcommand{\bea}{\begin{eqnarray}}
\newcommand{\eea}{\end{eqnarray}}
\newcommand{\bean}{\begin{eqnarray*}}
\newcommand{\eean}{\end{eqnarray*}}

\newcommand{\bdefin}{\begin{defin}}
\newcommand{\blemma}{\begin{lem}}
\newcommand{\bprop}{\begin{prop}}
\newcommand{\btheor}{\begin{thm}}
\newcommand{\bcoro}{\begin{cor}}
\newcommand{\bconj}{\begin{conj}}
\newcommand{\bdefprop}{\begin{defprop}}
\newcommand{\bexam}{\begin{example}}
\newcommand{\bitem}{\begin{leer}}
\newcommand{\eitem}{\end{leer}}
\newcommand{\edefin}{\end{defin}}
\newcommand{\elemma}{\end{lem}}
\newcommand{\eprop}{\end{prop}}
\newcommand{\etheor}{\end{thm}}
\newcommand{\ecoro}{\end{cor}}
\newcommand{\econj}{\end{conj}}
\newcommand{\brem}{\begin{rema}}
\newcommand{\erem}{\end{rema}}
\newcommand{\edefprop}{\end{defprop}}
\newcommand{\eexam}{\endexem\end{example}}

\def\1#1{{\bf #1}}
\def\2#1{{\cal #1}}
\def\3#1{{\sl #1}}
\def\4#1{{\tt #1}}
\def\5#1{{\sf #1}}
\def\6#1{{\mathfrak #1}}
\def\7#1{{\mathbb #1}}

\def\qed{\hfill{$\blacksquare$}\medskip}

\def\id{\mathrm{id}}

\newcommand{\End}{\mathrm{End}}
\newcommand{\Hom}{\mathrm{Hom}}

\newcommand{\Rep}{\mathrm{Rep}}
\newcommand{\Mod}{\mathrm{Mod}}

\newcommand{\obj}{\mathrm{Obj}}

\newcommand{\Aut}{\mathrm{Aut}}
\newcommand{\mcirc}{\,\circ\,}

\def\prf{\noindent \emph{Proof.\ }}

\newcommand{\restr}{\upharpoonright}
\newcommand{\rarr}{\rightarrow}

\newcommand{\impl}{\Rightarrow}
\newcommand{\ol}{\overline}

\newcommand{\ve}{\varepsilon}
\newcommand{\op}{{\mbox{\scriptsize op}}}
\newcommand{\rev}{{\mbox{\scriptsize rev}}}
\newcommand{\tr}{\mathrm{tr}}
\newcommand{\del}{\partial}

\newarrow{Congruent} 33333
\newarrow{Corresponds} <--->

\def\mobj#1{\raise .4\unitlens\hbox{\put(0,0){$#1$}}}

\numberwithin{equation}{section}

\begin{document}

\title{Modular Categories\footnote{To appear in: C. Heunen, M. Sadrzadeh, E. Grefenstette (eds.): {\it Compositional methods
in quantum physics and linguistics}, Oxford University Press, 2012.}}

\author{M.~M\"uger \\ IMAPP, Radboud University Nijmegen \\ The Netherlands}

\maketitle

\section{Introduction}
Modular categories, as well as the (possibly) more general non-degenerate braided fusion categories,
are braided tensor categories that are linear over a field and satisfy some natural additional
axioms, like existence of duals, semisimplicity, finiteness, and an important non-degeneracy condition. (Precise
definitions will be given later.) There are several reasons to study modular categories:
\begin{itemize}
\item As will hopefully become clear, they are rather interesting mathematical structures in
  themselves, well worth being studied for intrinsic reasons. For example, there are interesting
  number theoretic aspects.
\item Among the braided fusion categories, modular categories are the opposite extreme of the
  symmetric fusion categories, which are well known to be closely related to finite groups. Studying
  these two extreme cases is also helpful for understanding and classifying those braided fusion
  categories that are neither symmetric nor modular.
\item Modular categories serve as input datum for the Reshetikhin-Turaev construction of topological
  quantum field theories in $2+1$ dimensions and therefore give rise to invariants of smooth
  3-manifolds. This goes some way towards making Witten's interpretation of the Jones polynomial via
  Chern-Simons QFT rigorous. (But since there still is no complete rigorous non-perturbative
  construction of the Chern-Simons QFTs by conventional quantum field theory  methods, there also is
  no proof of their   equivalence to the   RT-TQFTs constructed using the representation theory of
  quantum groups.) 
\item Modular categories arise as representation categories of loop groups and, more generally, of
  rational chiral conformal quantum field theories. In chiral CQFT, the field theory itself, its
  representation category, and the conformal characters form a remarkably tightly connected
  structure. 
\item Also certain massive quantum field theories and quantum spin systems in two spatial
  dimensions lead to modular categories, e.g.\ Kitaev's `toric code'.
\item The recent topological approaches to quantum computing, while differing in details, all
  revolve around the notion of modular categories.
\end{itemize}

As the above list indicates, modular categories -- and related mathematical subjects like
representation theory of loop groups and of quantum groups at root-of-unity deformation parameter --
represent one of the most fruitful places of interaction of `pure' mathematics and mathematical
physics. 
(While modular categories play a certain r\^ole in string theory via their importance for rational
conformal field theories, the author believes that their appearance in massive field theories and
spin models may ultimately turn out to be of larger relevance for realistic physics.)

This article assumes a certain familiarity with category theory, including monoidal (=tensor)
categories. Concerning braided categories, only the definition, for which we refer to Majid's
contribution, will be assumed. Our standard reference for category theory is \cite{cwm}. For a
broader survey of some related matters concerning tensor categories, cf.\ also \cite{mue-arg}.


\section{Categories}
We limit ourselves to recalling some basic definitions.
A category is called an Ab-category if it is enriched over (the symmetric tensor category Ab of)
abelian groups. I.e., all hom-sets come with abelian group structures and the composition
$\circ$ of morphisms is a homomorphism w.r.t.\ both arguments. An Ab-category is called additive if
it has a zero object and every pair of objects has a direct sum. If $k$ is a field (or more
generally, a commutative unital ring) then a category $\2C$ is called $k$-linear if all hom-sets are
$k$-vector spaces (or $k$-modules) and $\circ$ is bilinear. An object $X$ is called simple if every
monic morphism $Y\hookrightarrow X$ is an isomorphism and absolutely simple if 
$\End\,X\cong k\id_X$. If $k$ is an algebraically closed field, as we will mostly assume, the two
notions coincide. A $k$-linear category is semisimple if every object is a finite direct sum of
simple objects. A semisimple category is called finite if the number of isomorphism classes of
simple objects is finite. (There is a notion of finiteness for non-semisimple categories,
cf.\ \cite{EO1}, but we will not need it.) A positive $*$-operation on a $\7C$-linear category $\2C$
is a contravariant endofunctor of $\2C$ that acts like the identity on objects, is involutive
($**=\id$) and anti-linear. (Dropping the $\7C$-linearity, one arrives at the notion of a
dagger-category.) A $*$-operation is called positive if $s^*\circ s=0$ implies $s=0$. A
category equipped with a (positive) $*$-operation is called hermitian (unitary). A unitary category
with finite-dimensional hom-spaces and splitting idempotents is semisimple (since finite dimensional
algebras with positive $*$-operation are semisimple, thus multi-matrix algebras).


\section{Tensor categories}
We assume familiarity with the basics of tensor (=monoidal) categories, including symmetric ones. We
therefore limit the discussion to (i) issues of duality in not necessarily braided 
tensor categories, (ii) fusion categories and (iii) module categories vs.\ categories of modules over an
algebra in a tensor category.

\subsection{Duality}
In the following, we will state definitions and results for strict tensor categories, but everything
can easily be adapted to the non-strict case.

\bdefin Let $\2C$ be a tensor category and let $X,Y\in\2C$. If there exist morphisms 
$e:Y\otimes X\rarr\11,\ d:\11\rarr X\otimes Y$ satisfying 
\[ (\id_X\otimes e)\circ(d\otimes\id_X)=\id_X, \quad\quad\quad
   (e\otimes\id_Y)\circ(\id_Y\otimes d)=\id_Y, \] 
then $Y$ is called a left dual of $X$ and $X$ is a right dual of $Y$. If every object has a left
(right) dual, we say that $\2C$ has left (right) duals.
\edefin

It is easy to prove that duals are unique up to isomorphism: If $Y,Y'$ are both left (or both right)
duals of $X$ then $Y\cong Y'$. This justifies writing ${}^\vee\!X$ ($X^\vee$) for a left (right)
dual of $X$. However, if $X$ admits a left dual ${}^\vee\!X$ and a right dual $X^\vee$, it may or
may not be the case that ${}^\vee\!X\cong X^\vee$. Clearly, a left dual ${}^\vee\!X$ of $X$ is also
a right dual of $X$ if and only if $X$ is a left dual of ${}^\vee\!X$, thus if and only if $X$ is
isomorphic to its second left dual ${}^{\vee\vee}\!X$. If ${}^\vee\!X\cong X^\vee$ holds, we say
that $X$ has a two-sided dual and mostly write $\ol{X}$ rather than ${}^\vee\!X$ or $X^\vee$. If all
objects have two-sided duals, we say that $\2C$ has {\it two-sided   duals}. 

There are three situations where duals, to the extent that they exist, are automatically two-sided: 
\begin{itemize}
\item[(i)]  $\2C$ is hermitian: If $({}^\vee\!X,e_X,d_X)$ defines a left dual of $X$, one finds that
$({}^\vee\!X,d_X^*,e_X^*)$ is a right dual. In this situation, some authors talk about conjugates
rather than duals and give slightly different axioms, cf.\ \cite{lro}.
\item[(ii)] $\2C$ is braided, cf.\ Section \ref{ss-db}.
\item[(iii)] $\2C$ is semisimple with simple unit, cf.\ \cite[Proposition 2.1]{eno}.
\end{itemize}

If $\2C$ is linear over a field $k$ and the unit $\11$ is absolutely simple,
we can and will use the bijection $k\rarr\End\,\11,\ c\mapsto
c\id_\11$ to identify $\End\,\11$ with $k$. If $X\in\2C$ is absolutely simple
and has a two-sided dual $\ol{X}$, one defines `squared dimension'
$d^2(X)=(e'\circ d)(e\circ d')$, which is easily seen to be
independent of the choices of the duality morphisms $e,d,e',d'$, cf.\ \cite{mue09}. 
If $X,Y$ and $X\otimes Y$ are absolutely simple, one has $d^2(X\otimes Y)=d^2(X)d^2(Y)$.  
If $\2C$ is semisimple and finite, we define $\dim\2C=\sum_i d^2(X_i)$, where $\{ X_i\}$ is a complete family of simple objects.

In the $k$-linear case, one would like to have a dimension function $X\mapsto d(X)$ that is additive and
multiplicative. If $\2C$ is semisimple
and finite, this can be done using Perron-Frobenius theory, cf.\ the next section. Another approach 
strengthens the requirement of existence of duals by introducing a new piece of structure:

\bdefin A left duality on a tensor category is an assignment $X\mapsto({}^\vee\!X,e_X,d_X)$, where
${}^\vee\!X$ is a left dual of $X$, the morphisms  
$e_X: {}^\vee\!X\otimes X\rarr\11,\ d_X:\11\rarr X\otimes{}^\vee\!X$ satisfying the above identities.
Similarly, a right duality is an assignment $X\mapsto(X^\vee,e'_X,d'_X)$ where the morphisms
$e'_X:X\otimes X^\vee\rarr\11,\ d'_X:\11\rarr X^\vee\otimes X$ establish $X^\vee$ as a right dual of $X$.
\edefin

If $\2C$ is a tensor category with a chosen left duality $X\mapsto({}^\vee\!X,e_X,d_X)$ and
$f:X\rarr Y$ then
\[ {}^\vee\!f=(e_Y\otimes\id_{{}^\vee\!X})\circ(\id_{{}^\vee\!Y}\otimes f\otimes\id_{{}^\vee\!X})
   \circ(\id_{{}^\vee\!X}\otimes d_X) : {}^\vee\!Y\rarr{}^\vee\!X \]
extends $X\mapsto{}^\vee\!X$ to a contravariant endofunctor ${}^\vee -$ of $\2C$.
Since ${}^\vee\! Y\otimes{}^\vee\! X$ is a left dual of $X\otimes Y$, ${}^\vee -$ can be considered as a
(covariant) tensor functor $\2C\rarr\2C^{\op,\rev}$, where $\2C^{\rev}$ coincides with $\2C$ as a category,
but has the  reversed tensor product $X\otimes^{\rev}Y=Y\otimes X$.
(Similarly for $-^\vee$.)
If ${}^\vee\!-$ is a left duality and $\2C$ has two-sided duals, one can find
a natural isomorphism $\gamma: \id\rarr {}^{\vee\vee}\!-$. Since this is of little use unless
$\gamma$ is monoidal, one defines:

\bdefin A pivotal category \cite{FY1,FY2} is a tensor category together with a left duality 
${}^\vee\!-$ and a monoidal natural isomorphism $\gamma: \id\rarr {}^{\vee\vee}\!-$.
\edefin

\brem 1. Categories equipped with left and/or right dualities are often called rigid or autonomous, and
pivotal categories are also called sovereign. We will avoid all these terms. The categories we
consider will either have the property of possessing two-sided duals (without given duality
structures) or be equipped with pivotal or spherical structures.

2. For a general tensor category with two-sided duals, there is little {\it a priori} reason to expect the existence
of a pivotal structure, but see Theorems \ref{theor-eno4} and \ref{theor-eno} below.
\erem

\bdefin Let $\2C$ be a strict pivotal category with left duality $X\mapsto(\ol{X},e_X,d_X)$ and
monoidal natural isomorphism $(\alpha_X:X\rarr\ol{\ol{X}})_{X\in\2C}$. For $X\in\2C$ and
$s\in\End\,X$, define the left and right traces of $s$ by
\bean \mbox{tr}^L_X(s) &=& e_X\circ(\id_{\ol{X}}\otimes (s\circ\alpha^{-1}_X))\circ d_{\ol{X}},\\
  \mbox{tr}^R_X(s) &=& e_{\ol{X}}\circ((\alpha_X\circ s)\otimes\id_{\ol{X}})\circ d_X.
\eean
(Both traces take values in $\End\,\11$.) The left and right dimensions of $X$ are defined by $d_{L/R}(X)=tr^{L/R}_X(\id_X)$.
\edefin

The traces satisfy $\tr_X(s\circ t)=\tr_Y(t\circ s)$ for all $t:X\rarr Y$ and $s:Y\rarr X$, as well as
$\tr_{X\otimes Y}(s\otimes t)=\tr_X(s)\tr_Y(t)$ for $s\in\End\,X,\ t\in\End\,Y$. If $X$ is
absolutely simple, the dimensions defined in terms of the traces are connected to the intrinsic squared
dimension by $d_L(X)d_R(X)=d^2(X)$.

A pivotal category is called spherical \cite{BW1} if $\tr_X^L=\tr_X^R$ for all $X\in\2C$. In this
case we have $d_L(X)=d_R(X)$ and simply write $d(X)$. In fact, a pivotal category that is semisimple
over a field is
spherical if and only if $d_L(X)=d_R(X)$ for every $X$. To this day, all known examples of pivotal
fusion categories are spherical. (E.g., if $H$ is a finite dimensional semisimple Hopf algebra in
characteristic zero, one automatically has $S^2=\id$ and therefore sphericity of the module
category.) 

As noted before, in a hermitian category, duals (if they exist) are always two-sided. In a unitary
category there is a canonical way of defining traces of endomorphisms, giving rise to positive
dimensions, that makes no use of spherical structures, cf.\ \cite{lro}. But as shown by Yamagami
\cite{yamag}, a unitary 
category always admits a unique spherical structure that gives rise to these traces. 
Note, however, there are $\7C$-linear categories that do not admit a unitary structure,
cf. \cite{row1}. (Essentially, quantum group categories at odd root of unity.)

\subsection{Fusion categories}\label{ss-fusion}
Let $k$ be an algebraically closed field. A $k$-linear tensor category is called fusion category if
it has finite dimensional hom-spaces, is semisimple with finitely many isomorphism classes of simple
objects, the unit $\11$ is absolutely simple and all objects have duals (which are automatically
two-sided by semisimplicity).
A fusion subcategory of a fusion category is a full tensor subcategory that is again
fusion. (I.e.\ closed under direct sums and duals.)

Even for fusion categories it is unknown whether pivotal structures always exist. (See \cite{DSSP}
for work currently in progress.) But there is a result in this direction:

\btheor \label{theor-eno4} \cite{eno}
Let $\2C$ be a fusion category over $\7C$ and let ${}^\vee -$ be an arbitrary left dual structure. Then
the tensor functor ${}^{\vee\vee\vee\vee}\!-$ is naturally
monoidally isomorphic to the identity functor.
\etheor

In a fusion category, we can find mutually non-isomorphic simple objects $X_i$, indexed by a finite
set $I$, such that every object is a finite direct sum of the objects $X_i, i\in I$. There is a
distinguished index $0\in I$ such that $X_0=\11$. Now we can define non-negative integers
$N_{ij}^k\in\7N_0$ via $X_i\otimes X_j\cong\oplus_{k\in I} N_{ij}^k X_k$. These numbers have various
obvious properties, like $N_{i0}^j=N_{0i}^j=\delta_{i,j}$. We also have
$N_{ij}^0=\delta_{i,\ol{\jmath}}$, where $i\mapsto\ol{\imath}$ is the involution on $I$ defined by
$\ol{X_i}\cong X_{\ol{\imath}}$. The structure $(I,0,i\mapsto\ol{\imath}, N_{\cdot,\cdot}^\cdot)$
is known as a discrete hypergroup. (In particular, every group $G$ gives rise to a hypergroup,
taking $I=G,\, 0=e,\, \ol{g}=g^{-1}$ and $N_{g,h}^k=\delta_{gh,k}$.)

For every discrete hypergroup $(I,\cdots)$ there
is a unique group $G(I)$ equipped with a map $\del:I\rarr G(I)$ satisfying 
\[ \del 0=e,\ \  \del\ol{\imath}=(\del i)^{-1}, \quad\mbox{and}\quad N_{i,j}^k>0\impl \del k=\del i\cdot \del j\]
and being universal for such maps. (Thus for every map $\del':I\rarr H$, where $H$ is a group, satisfying the
same axioms as $\del$ there is a unique group homomorphism $\alpha:G\rarr H$ such that
$\del'=\alpha\circ\del$.) In view of this property, $G(I)$ is called the (universal) grading group
of the hypergroup $I$. (Also `groupification' might be a good term, in analogy to abelianization.)
Cf.\ \cite{BL,GN1}. Notice that $G(I)$ is abelian when $\2C$ is braided. If $H$ is a compact group,
let $\widehat{H}$ be the hypergroup corresponding to the semisimple category $\Rep\,H$. Then
the (discrete abelian) grading group $G(\widehat{H})$ is canonically isomorphic to Pontrjagin dual
of the (compact abelian) center $Z(H)$ of $H$, cf.\ \cite{mue14}.

Each matrix $N_i=(N_{ij}^k)_{j,k}$ is irreducible and has
non-negative entries. Thus it has a unique positive Perron-Frobenius eigenvalue which is denoted as
FP$\dim(X_i)$, the Frobenius-Perron dimension of $X_i$, cf.\ \cite{eno}. The Frobenius-Perron
dimension \cite{eno} of the fusion category $\2C$ is defined by 
\[ \mathrm{FP}\dim(\2C)=\sum_{i\in I} \mathrm{FP}\dim(X_i)^2. \]

By definition, FP$\dim(X)$ and FP$\dim(\2C)$ live in $\7R_{>0}$ rather than in the ground
field $k$, but if $\7R\subset k$ it makes sense to compare FP$\dim(X)$ with $d^2(X)$ (which is
canonically defined for all simple objects) and with $d(X)$ (if $\2C$ has a spherical structure).
It is classical that the FP-dimension is the unique positive dimension function 
on a finite hypergroup $I$. In particular the Frobenius-Perron dimension coincides with the positive
dimension function defined on unitary categories with duals \cite{lro}. (Recall that the latter
arises from a unique spherical structure.) Somewhat more generally, if $\2C$ is spherical over $k\supset\7R$
and $d(X)>0$ for all $X$ then $d(X)=\mathrm{FP}\dim(X)$ for all $X$. 
Since not every fusion category over $\7C$ is unitarizable (in the sense of admitting a positive
$*$-operation) it is important that there are the following remarkable results:

\btheor \label{theor-eno}
\cite{eno} Let $\2C$ be a fusion category over $\7C$. Then
\begin{itemize}
\item[(i)] For every simple $X$, one has $0<d^2(X)\le\mathrm{FP}\dim(X)^2$. Thus
  $1\le\dim\2C\le\mathrm{FP}\dim\2C$.
\item[(ii)] If $\dim\2C=\mathrm{FP}\dim(\2C)$ (equivalent to $d^2(X)=\mathrm{FP}\dim(X)^2$ for every simple
  $X$) $\2C$ admits a unique spherical structure for which $d(X)=\mathrm{FP}\dim(X)>0$ for every
  $X$. Such categories are called pseudo-unitary.
\end{itemize}
\etheor

\subsection{Algebras in tensor categories and their modules}\label{ss-algebras1}
A considerable part of (commutative) algebra can be generalized from the symmetric categories Ab and
Vect$_k$ to arbitrary (braided) tensor categories. This plays an important r\^ole in the structural
study of such categories and in particular of braided fusion categories. In this subsection we
discuss some facts that do not require a braiding.

\bdefin Let $\2C$ be a strict tensor category. An algebra (or monoid) in $\2C$ is a triple
$(A,m,\eta)$, where $A\in\2C$, $m:A\otimes A\rarr A$ and $\eta:\11\rarr A$ satisfy 
$m\circ(m\otimes\id_A)=m\circ(\id_A\otimes m)$ (associativity) and
$m\circ(\eta\otimes\id_A)=\id_A=m\circ(\id_A\otimes\eta)$ (unit property). If $\2C$ is non-strict,
one inserts associativity isomorphisms at the appropriate places.
\edefin

\brem At least in categories that are not linear over a field, it would be more appropriate to speak
of monoids rather than algebras, but the latter term is used much more in the recent literature.
\erem

\bdefin If $\2C$ is a tensor category and $(A,m,\eta)$ is an algebra in $\2C$ then a left $A$-module
(in $\2C$) is a pair $(X,\mu)$ where $X\in\2C$ and $\mu:A\otimes X\rarr X$ satisfies 
$\mu\circ(m\otimes\id_X)=\mu\circ(\id_A\otimes\mu)$ and $\mu\circ(\eta\otimes\id_X)=\id_X$.
\[ \begin{tangle} \step[1.5]\object{X}\\ \hstep\ld\step[.3]\obj{\mu}\\ \hh\step[-.7]\obj{m}\step[.7]\cd\hstep\id\\ 
   \object{A}\step\object{A}\hstep\object{X} \end{tangle}\quad=\quad
   \begin{tangle} \step[2]\object{X}\\ \Ld\step[.3]\obj{\mu}\\ \id\step\ld\step[.3]\obj{\mu}\\  
\object{A}\step\object{A}\step\object{X} \end{tangle}
\quad\quad\quad\quad\quad\quad
\begin{tangle} \step\object{X}\\ \ld\step[.3]\obj{\mu}\\ \obj{\eta}\counit\step\id\\ \step\object{X} \end{tangle}
   \quad\quad=\quad\quad   \begin{tangle}\object{X}\\\id\\\object{X}\end{tangle}
\]

The left $A$-modules in $\2C$ form a category, denoted $A$-$\Mod_\2C$ or ${}_A\2C$, with hom-sets
\[ \Hom_{{}_A\2C}((X,\mu),(X',\mu'))=\{ t\in\Hom_\2C(X,X')\ | \ t\circ\mu=\mu'\circ(\id_A\otimes t)\}.\]
Right modules and bimodules are defined analogously. 
\edefin

There is a functor $F_A:\2C\rarr {}_A\2C$, the free module functor, defined by 
$F_A: X\mapsto (A\otimes X, m\otimes\id_X)$. Notice that for every $(X,\mu)\in{}_A\2C$ we have
$\mu\in\Hom_{{}_A\2C}(F(X),(X,\mu))$. The functor $F_A$ is faithful provided $s\mapsto\id_A\otimes s$ is
injective, which usually is the case. Since the maps
\bean \Hom_{{}_A\2C}(F_A(X),F_A(Y)) \rarr \Hom_\2C(A\otimes X,Y), && s\mapsto s\circ(\eta\otimes\id_Y), \\
   \Hom_\2C(A\otimes X,Y)\rarr\Hom_{{}_A\2C}(F_A(X),F_A(Y)), && t\mapsto (m\otimes\id_Y)\circ(\id_A\otimes t)\eean
are inverses of each other, we have a bijection $\Hom_{{}_A\2C}(F_A(X),F_A(Y)) \cong \Hom_\2C(A\otimes X,Y)$.
Thus $F_A$ will in general not be full, and it can happen that $F_A$ trivializes an object $X\in\2C$
in the sense of mapping it to a multiple of the unit object $(A,m)$ of ${}_A\2C$.

The free $A$-modules form a full subcategory of ${}_A\2C$, but in order to say more about the module
category, one needs a descent-type assumption like the following:

\bdefin An algebra $(A,m,\eta)$ in a tensor category is called separable if the
multiplication morphism admits a splitting that is a morphism of $A$-$A$ bimodules, i.e.\ a morphism
$\widetilde{m}: A\rarr A\otimes A$ 
satisfying 
\[ (m\otimes\id_A)\circ(\id_A\otimes\widetilde{m})=\widetilde{m}\circ m
   =(\id_A\otimes m)\circ (\widetilde{m}\otimes\id_A), \quad\quad\quad
 m\circ\widetilde{m}=\id_A. \]
\[ \begin{tangle} \hh\hcd\step\id \\ \hh\id\step\hcu\end{tangle}\quad=\quad
   \begin{tangle}\hh\hcu\\ \hh\hcd\end{tangle}\quad=\quad
   \begin{tangle} \hh\id\step\hcd \\ \hh\hcu\step\id\end{tangle}
\quad\quad\quad\quad\quad\quad\quad
  \begin{tangle} \hh\hcd\\ \hh \hcu\end{tangle}\quad=\quad\begin{tangle}\id\end{tangle}
\]
If $\2C$ is $k$-linear and $\dim\Hom(\11,A)=1$ then the algebra $(A,m,\eta)$ is called connected.
\edefin

\blemma If $\2C$ is a tensor category and $(A,m,\eta)$ a separable algebra, then every module
$(X,\mu)\in{}_A\2C$ is a quotient of the free module $F(X)$.
\elemma 

\prf Let $\widetilde{m}:A\rarr A\otimes A$ be a splitting of $m$. Defining
\[ \gamma=(\id_A\otimes\mu)\circ(\widetilde{m}\otimes\id_X)\circ(\eta\otimes\id_X): X\rarr A\otimes X, \]
an easy computation shows that $\gamma\in\Hom_{{}_A\2C}((X,\mu),F(X))$ and
$\mu\circ\gamma=\id_X$. Thus $\mu\in\Hom_{{}_A\2C}(F(X),(X,\mu))$ is a split epimorphism.
\qed

\brem A notion similar to separability appeared in \cite{brug}, where a morphism 
$\beta:\11\rarr A\otimes A$ was required, satisfying axioms following from the above ones if one
takes $\beta=\widetilde{m}\circ\eta$. (Notice that $\beta$ is what is actually used in the Lemma.)
A related concept is that of a special Frobenius algebra, which is a quintuple
$(A,m,\eta,\widetilde{m},\ve)$ where $(A,m,\eta,\widetilde{m})$ is a separable algebra and
$(A,\widetilde{m},\ve)$ is a coalgebra. (In particular, $\widetilde{m}$ must be coassociative.) 
\erem

\blemma \cite{ostrik,eno,dmno} Let $\2C$ be a fusion category and $(A,m,\eta)$ an algebra in $\2C$. Then the
following are equivalent:
\begin{itemize}
\item[(i)] $A$ is separable.
\item[(ii)] The category ${}_A\2C$ is semisimple.
\item[(iii)] The category $\2C_A$ is semisimple.
\item[(iv)] The category ${}_A\2C_A$ is semisimple.
\end{itemize}
\elemma

\brem For $\2C$ spherical, (i)$\impl$(iv) was already shown in \cite{mue09}, as part of the following. \erem

\btheor \cite{mue09} If $\2C$ is a spherical fusion category and $A\in\2C$ a separable and connected
algebra then $\2D={}_A\2C_A$ is a spherical fusion category and it contains a separable connected algebra $B$ such that
$\2C\simeq{}_B\2D_B$. Calling fusion categories that are related in this way {\it weakly monoidally
  Morita equivalent}, weak monoidal Morita equivalence is an equivalence relation, denoted by $\approx$.
If $\2C\approx\2D$ then $\dim\2C=\dim\2D$.
\etheor

\brem 1. If $H$ is a finite dimensional semisimple and co-semisimple Hopf algebra then
$H-\Mod\approx\widehat{H}-\Mod$.

2. The results of the theorem have been generalized to fusion categories in \cite{eno}. In that
generality, one must use FP$\dim$ instead of $\dim$.
\erem

\subsection{Module categories and categories of modules}
\bdefin \label{def-Z0}
If $\2M$ is a category, $\End\,\2M$ denotes the category whose objects are the endofunctors
of $\2M$, i.e.\ functors $\2M\rarr\2M$, and whose morphisms are the natural transformations of these
endofunctors. With $\otimes$ defined as composition of functors, $\End\,\2M$ is a strict tensor
category. An alternative name for $\End\,\2M$ is $Z_0(\2M)$, the monoidal center of $\2M$. 
\edefin 

\bdefin Let $\2C$ be a tensor category and $\2M$ a category. A left, resp.\ right, $\2C$-module structure on
$\2M$ is a tensor functor $\2C\rarr\End\,\2M$, resp.\ $\2C^{\rev}\rarr\End\,\2M$. A left (right) $\2C$-module category is a
category $\2M$ together with a left (right) $\2C$-module structure. A module category is
indecomposable if it is not the direct sum of two non-trivial module categories.

Left (right) $\2C$-module categories form a 2-category, whose 1-morphisms are functors between
module categories intertwining the $\2C$-actions and whose 2-morphisms are natural transformations
between module functors.
\edefin

\brem Unpacking the above definition, one finds that an action of a tensor category $\2C$ on a
category $\2M$ is the same as a functor $\2C\times\2M\rarr\2M$ satisfying certain properties like
the existence of natural isomorphisms $(C\otimes D)\otimes M\rarr C\otimes(D\otimes M)$, etc. For
such an approach, cf.\ e.g.\ \cite{ostrik}.
\erem

If $\2C$ is a tensor category and $(A,m,\eta)$ an algebra in $\2C$ then the left module category
${}_A\2C$ has a structure as a right $\2C$-module category: $(X,\mu)\otimes Y=(X\otimes
Y,\mu\otimes\id_Y)$. Similarly, $\2C_A$ has a left $\2C$-module structure. The question now arises
whether every $\2C$-module category is of this form. For fusion categories this is the case:

\btheor \cite{ostrik}
Let $\2C$ be a fusion category and $\2M$ an indecomposable semisimple left module category over
$\2C$. Then there is a connected separable algebra $(A,m,\eta)$ in $\2C$ and an equivalence
$\2M\rarr\2C_A$ of module categories. 
\etheor

Thus (indecomposable, semisimple) module categories over a fusion category $\2C$ and categories of
modules for a separable algebra in $\2C$ are essentially the same thing. This allows to discuss the
weak monoidal Morita equivalence of \cite{mue09} in terms of module categories: If $\2M$ is a good
module category over the fusion category $\2C$ then the category $\Hom_\2C(\2M,\2M)$ of $\2C$-module
functors, which clearly is monoidal, actually is fusion and is called dual to $\2C$ with respect to
$\2M$, cf.\ \cite{eno}.  The connection with \cite{mue09} is as follows: When $\2M={}_A\2C$ then
$\Hom_\2C(\2M,\2M)\simeq {}_A\2C_A$, where ${}_A\2C_A$ acts on ${}_A\2C$ by tensoring from the left.


\section{Braided tensor categories}
\subsection{Centralizers in braided categories and the symmetric center $Z_2$}
In this section, we quickly discuss some aspects of braided tensor categories that do not require
additional axioms like existence of duals or linearity over a ring or field. A braided tensor category is a
tensor category $\2C$ equipped with a braiding $c$. For the definition of the latter, we refer to
Majid's contribution to the present volume.
In principle, a braided tensor category should be written as $(\2C,c)$, where $\2C$ is a tensor
category and $c$ a braiding on $\2C$, but we will suppress the $c$. (After all, one usually does the
same with the various items of the monoidal structure.)

\bdefin \label{def-Ctil}
Let $\2C$ be a braided tensor category with braiding $c$. The opposite braiding
$\widetilde{c}$ is defined by $\widetilde{c}_{X,Y}=(c_{Y,X})^{-1}$. 
The tensor category $\2C$ equipped with the braiding $\widetilde{c}$ is denoted $\widetilde{\2C}$.

In the graphical calculus (where we draw morphisms going upwards) the braidings $c$ and
$\widetilde{c}$ are represented by 
\[ c_{X,Y}=\begin{tangle} \object{Y}\Step\object{X}\\ \xx\\ \object{X}\Step\object{Y} \end{tangle}
 \quad\quad \quad\quad
   \widetilde{c}_{X,Y}=\begin{tangle} \object{Y}\Step\object{X}\\ \x\\ \object{X}\Step\object{Y} \end{tangle}\] 
which is consistent in view of the isotopy
\[ c_{Y,X}\circ \widetilde{c}_{X,Y}=\quad \begin{tangle}
  \xx\\ \x \end{tangle}\quad=\quad \begin{tangle}\id\step\id\end{tangle}\quad=\id_{X\otimes Y}. \]
\edefin

\brem 1. Clearly $\widetilde{\widetilde{c}}=c$ always holds. The stronger statement
$\widetilde{c}=c$ is equivalent to
\be c_{X,Y}\circ c_{Y,X}=\id_{Y\otimes X} \ \ \forall X,Y\in\2C. \label{eq-sym}\ee
Recall that a BTC satisfying this additional condition is called a symmetric tensor category (STC).
For this reason the braiding (=symmetry) of an STC is depicted by a crossing with unbroken
lines. Conversely, the definition of BTCs is obtained from the usual definition of symmetric tensor categories by
dropping the condition (\ref{eq-sym}) that the braiding be involutive. (In doing so, one must add a
second hexagon axiom which follows from the first in the presence of (\ref{eq-sym}).)

2. Since examples of STCs abound, e.g.\ the category of sets, categories of modules and vector spaces,
representation categories of groups, categories of sheaves, etc., it is quite astonishing that they
were formalized only in 1963, cf.\ \cite{benab,macl2}. Considering that it is much harder to find
interesting examples of non-symmetric BTCs, it is less surprising that their formalization took place 
only around 1985/6 in the first preprint versions of \cite{js2}. To some extent, this was
inspired by the surge of activity around new `quantum' invariants in low dimensional topology
\cite{jones,HOMFLY}, quantum groups \cite{drin}, conformal field theory \cite{ms1} and algebraic quantum
field theory \cite{frs} in the second half of the 1980s. But there are also considerations
intrinsic to (higher) category theory as well as in `old-fashioned' (non-`quantum') algebraic
topology that lead to BTCs, cf.\ \cite{js2}. 

3.  In a BTC, any two objects commute up to isomorphism: $X\otimes Y\cong Y\otimes X$. 
Looking for a stronger statement to be attached to the expression `$X$ and $Y$ commute', one could
think of $c_{X.Y}=\id_{X\otimes Y}$ or the slightly weaker $X\otimes Y=Y\otimes X$. But these are
hardly ever satisfied in interesting BTCs. This essentially leaves us with only the following option.
\erem

\bdefin Let $\2C$ be a BTC.
\begin{itemize}
\item[(a)] Two objects $X,Y\in\2C$ are said to commute if $c_{X,Y}=\widetilde{c}_{X,Y}$,
  equivalently $c_{X,Y}\circ c_{Y,X}=\id_{Y\otimes    X}$. 
\item[(b)] Let $\2D\subset\2C$ be a subcategory. The centralizer $C_\2C(\2D)$ is the full
subcategory of $\2C$ defined by
\[ \obj\, C_\2C(\2D)= \{ X\in\2C\ | \ c_{X,Y}\circ c_{Y,X}=\id_{Y\otimes X} \ \forall Y\in\2D \}.\]
\item[(c)] The symmetric center $Z_2(\2C)$ is defined as $C_\2C(\2C)$, i.e.\ the full subcategory
  defined by
\[ \obj\, Z_2(\2C)= \{ X\in\2C\ | \ c_{X,Y}\circ c_{Y,X}=\id_{Y\otimes X} \ \forall Y\in\2C \}.\]
\end{itemize}
\edefin

\brem 1. The definition of $Z_2(\2C)$ is due to \cite{brug,mue06}, but the concept first appeared
ten years earlier in the context of algebraic quantum field theory \cite{khr}. The centralizer
$C_\2C(\2D)$ seems first to have appeared in \cite{mue11}.

2. In most of the literature, $\2D'$ is written instead of $C_\2C(\2D)$. When the ambient
category $\2C$ is fixed, there is no risk of confusion.

3. $C_\2C(\2D)$ depends only on the objects of $\2D$. One easily sees, for any $\2D$, that the
monoidal unit $\11$ lies in $C_\2C(\2D)$ and that 
$C_\2C(\2D)$ is closed under tensor products, thus is is a full monoidal subcategory. Furthermore,
it is closed under isomorphisms (i.e.\ replete) and under direct sums, if they exist. The braiding
that $C_\2C(\2D)$ inherits from $\2C$ in fact is a symmetry when $\2D=\2C$, thus $Z_2(\2C)$ is an STC.
In fact, a BTC $\2C$ is symmetric if and only if $\2C=Z_2(\2C)$. 

4. The objects of $Z_2(\2C)$ are called {\it  transparent} \cite{brug}, since for braidings
involving them it there is no difference between over- and under-crossings, or {\it central}.

5. The name `center' for $Z_2(\2C)$ is amply justified by considerations from
higher category theory, cf.\ e.g.\ \cite{baez,baezM}, some of which already played a r\^ole in
\cite{js2}. In higher category theory, an infinite family of center constructions is considered. In
the more limited context of 1-categories, one deals with the bicategories consisting of
categories, tensor categories, braided tensor categories, symmetric tensor categories,
respectively. Notice that moving rightwards in this list adds a piece of structure (tensor
structure, braiding) except in the last step, where a condition is added, to wit (\ref{eq-sym}).
It is clear that there are forgetful (2-)functors moving leftwards in this list. 
More interestingly, there are constructions, called centers, in the opposite direction. 
We have just defined the symmetric center $Z_2(\2C)$ of a BTC $\2C$, and the monoidal
center $Z_0(\2C)=\End\,\2C$ of a category $\2C$ was given in Definition \ref{def-Z0}. The braided
center $Z_1(\2C)$ of a tensor category $\2C$ will be defined below. (These constructions are 
categorifications of their simpler analogues for 0-categories, i.e.\ sets, where we deal with the
1-categories of sets, monoids and commutative monoids. The center $Z_1(M)$ of a monoid is well
known, whereas the center $Z_0(S)$ of a set $S$ is the monoid of endomaps of $S$.)

Notice that the center constructions are compatible with equivalences of categories, but not with
more general functors. Thus they are functorial only on the sub-bicategories of CAT, $\otimes$-CAT,
etc., whose 1-morphisms are equivalences of categories, tensor categories (TCs), etc.

6. In view of the above, one realizes that symmetric tensor categories, which play a rather prominent
r\^ole in large parts of mathematics, are but one extreme case of braided tensor categories, singled
out by the condition that they coincide with their symmetric centers. This makes it natural ask
whether interesting things can be said in the opposite extreme case, namely when the center of a BTC
$\2C$  is trivial in the sense of containing only what it must contain, to wit the unit object and its
direct sums. (Compare with the theory of von Neumann algebras, where the commutative ones and those
with trivial center (`factors'), play distinguished r\^oles.)
This is indeed the case, modular categories just being braided categories having a bit
more structure and having trivial symmetric center $Z_2(\2C)$. Since this section is devoted to
results requiring no additional axioms, the study of modular categories will begin later.
\erem

\subsection{Rambling remarks on the construction of proper BTCs}
So far, we have not given any example of a non-symmetric BTC. The simplest one, the free braided
tensor category $\2B$ generated by one object, or just the braid category, is constructed from the
braid groups $B_n$ and could have been found long before \cite{js2}. Its objects are the
non-negative integers $\{0, 1, 2, \ldots\}$ with addition as tensor product. The category is discrete,
i.e.\ $\Hom_\2B(n,m)=\emptyset$ when $n\ne m$, the endomorphisms given by $\End_\2B(n)=B_n$ with
composition as in $B_n$ (i.e.\ concatenation of braids). The tensor product of morphisms is given by
horizontal juxtaposition of braids, and the braiding is defined as in
\[ c_{3,2}=  \quad
\begin{picture}(100,50)(0,30)
\put(0,0){\line(1,1){70}}
\put(30,0){\line(1,1){70}}
\put(60,0){\line(1,1){70}}

\put(85,0){\line(-1,1){10}}
\put(70,15){\line(-1,1){10}}
\put(55,30){\line(-1,1){10}}
\put(40,45){\line(-1,1){25}}

\put(115,00){\line(-1,1){20}}
\put(85,30){\line(-1,1){10}}
\put(70,45){\line(-1,1){10}}
\put(55,60){\line(-1,1){10}}
\end{picture}
\]\\ \\ \\

This category clearly is not symmetric. (For example, under the isomorphism $B_2\rarr\7Z$, the braid 
$c_{1,1}\circ c_{1,1}$ is mapped to $2\ne 0$.)

An attempt to systematize the known constructions of braided categories was made in \cite{mue-arg},
where three types of constructions were distinguished: 
\begin{enumerate}
\item Braided deformations of symmetric tensor categories.
\item Free (=topological) constructions.
\item The braided center $Z_1(\2C)$ of a tensor category.
\end{enumerate}
While the philosophies behind these three approaches are quite different, they are by no means
mutually exclusive. In fact, the most interesting braided categories, to wit the representation
categories of quantum groups, can be understood in terms of all three constructions! 

Space constraints do not allow to say much about the deformation
approach. While it is usually formulated in terms of a `$q$-deformation' of the universal enveloping
algebra $U(\6g)$ of a simple Lie algebra, giving rise to a quasi-triangular Hopf algebra $U_q(\6g)$,
cf.\ e.g.\ \cite{majid, kassel, cp}, one may argue that (as always?) a categorical perspective
provides additional insight. Namely, the representation categories $\2C(\6g,q)=U_q(\6g)-\Mod$ can be  
obtained directly by deformation of the STCs $\2C(\6g)=U(\6g)-\Mod$. Such deformations are
controlled by the third Davydov-Yetter cohomology, cf.\ \cite{dav0,yetter4,yetter}, and one can
show \cite{eno} that $H^3(D(\6g)-\Mod)$ is one-dimensional for a simple Lie algebra $\6g$, explaining
the one-parameter family of $q$-deformations. Actually constructing the braided deformed category
$\2C(\6g,q)$ from $\2C(\6g)$ can be done formally using Drinfeld associators \cite{kt2} or analytically
(i.e.\ non-formally) using the Knizhnik-Zamolodchikov connection \cite{kl}. The categories thus
obtained can be shown to be equivalent to the representation categories of the corresponding quantum
groups. 

The approach via `free constructions' generalizes the construction of the braid category $\2B$ given
above. Given a tensor category $\2C$, there is a free braided tensor category $F\2C$ over $\2C$,
which reduces to $\2B$ when $\2C$ is the trivial tensor category $\{\11\}$, cf.\ \cite{js2}. This
construction provides a left adjoint to the forgetful 2-functor from the bicategory of braided
tensor categories to the bicategory of tensor categories. There are analogous versions of this
construction in case one studies categories with additional structures, like duals or linearity over
a field. Most important, at least as far as connections with low dimensional topology are concerned,
are the the various categories of tangles, which can be considered as free rigid braided category,
free ribbon category, etc., generated by one object, cf.\ \cite{turaev0,FY1,turaev,kassel,yetter}. 
The tangle categories are not linear over a field, but can easily be linearized using the
the free vector space functor from sets to vector spaces. The categories thus obtained are still too
generic and too big (in the sense of having infinite dimensional hom-spaces) to be really
interesting, but quotienting them by a suitable ideal defined, e.g., in terms of a link invariant,
one can obtain rigid braided or ribbon categories with finite dimensional hom-spaces. For appropriate choices
of the link invariant (HOMFLY or Kauffman polynomials), one actually obtains the representation
categories of the quantum groups of types A-D, cf.\ \cite{tw2,bebl1,bebl2,blanchet}, including the
most interesting (and difficult) root-of-unity case. Again, there is no space to go into this any further.

However, beginning in the next subsection, we will have more to say about the third approach to the
construction of modular categories, the braided center construction, since it is of considerable
relevance for the structure theory of modular categories. This construction will play a central
r\^ole in most of what follows. 

\subsection{The braided center $Z_1(\2C)$} \label{ss-Z1}
For simplicity, we give the following definition only for strict tensor categories, but the
generalization is straightforward.

\bdefprop Let $\2C$ be a strict tensor category.
\begin{itemize}
\item[(a)] Let $X\in\2C$. A half braiding $e_X$ for $X$ is a family 
$\{e_X(Y): X\otimes Y\stackrel{\cong}{\longrightarrow}Y\otimes X\}_{Y\in\2C}$ 
of isomorphisms, natural w.r.t.\ $Y$, satisfying  $e_X(\11)=\id_X$ and
\[ e_X(Y\otimes Z)=\id_{Y}\otimes e_X(Z) \mcirc e_X(Y)\otimes\id_{Z} \quad \forall Y,Z\in\2C. \]
\item[(b)] Let $Z_1(\2C)$ be the category whose objects are pairs $(X,e_X)$ consisting of an object
  and a half-braiding, the hom-sets being given by
\[ \Hom_{Z_1(\2C)}((X, e_X),(Y, e_Y))= \{ t\in\Hom_\2C(X,Y) \ | \  \id_X\otimes t \mcirc e_X(Z)=
   e_Y(Z)\mcirc t\otimes\id_X \quad \forall Z\in\2C \}.\] 
Now a tensor product of two objects is defined by
$(X, e_X)\otimes(Y, e_Y)=(X\otimes Y, e_{X\otimes Y})$, where
\[ e_{X\otimes Y}(Z)= e_X(Z)\otimes\id_Y\mcirc\id_X\otimes e_Y(Z). \]
The tensor unit is $(\11,e_\11)$ where $e_\11(X)=\id_X$. Defining composition and tensor
product of morphisms to be inherited from $\2C$, one verifies that $Z_1(\2C)$ is a strict tensor
category. Finally,
\[  c_{(X, e_X),(Y, e_Y)}=  e_X(Y) \]
defines a braiding. The braided tensor category $Z_1(\2C)$ will be called the braided center of
$\2C$. 
\end{itemize}
\edefprop

\brem \label{rem-Z1} 
1. Usually, $Z_1(\2C)$ is denoted by $Z(\2C)$. We wrote $Z_1$ to avoid confusion with 
the symmetric center $Z_2$ (if $\2C$ is braided, both $Z_1(\2C)$ and $Z_2(\2C)$ are defined) but
later we will drop the subscript and identify $Z=Z_1$. 

2. In the same way as $Z_2(\2C)=C_\2C(\2C)$ is a special case of the centralizer $C_\2C(\2D)$, there
is a `relative center' $Z_1(\2C,\2D)$ for a pair $\2D\subset\2C$. Its objects are pairs $(X,e_X)$
where $X\in\2C$ and $e_X$ is a family of isomorphisms $X\otimes Y\rarr Y\otimes X$ for all
$Y\in\2D$. $Z_1(\2C,\2D)$ is always monoidal, but not necessarily braided, and we have
$Z_1(\2C,\2C)=Z_1(\2C)$. 

3. The definition of $Z_1$ appeared in \cite{js1} and \cite{majid1}. The second reference also
gave $Z(\2C,\2D)$ and attributed $Z_1(\2C)$ to unpublished work of Drinfeld, which led many authors
to call $Z_1(\2C)$ the `Drinfeld center'.

4. Despite its being somewhat involved, the definition of $Z_1(\2C)$ is quite natural. It is
nevertheless instructive to give an interpretation in terms of bicategories. The point is that a
(strict) tensor category is `the same' as a (strict) 2-category with one 
object, and similarly, braided tensor categories correspond to monoidal 2-categories with one
object \cite{js2}. Now, let $\2E$ be the 2-category with one object corresponding to the tensor
category $\2C$ and let $\2F$ be the monoidal 2-category $Z_0(\2E)$ of endo-2-functors of $\2E$. If
$\2F_1\subset\2F$ is the full sub-2-category retaining only the unit object $\11=\id_\2E$, it turns
out that the braided category corresponding to $\2F_1$ is nothing but $Z_1(\2C)$.

5. The construction of $Z_1(\2C)$ was preceded and probably motivated by Drinfeld's definition of
the quantum double $D(H)$, which is a quasi-triangular Hopf algebra, of a Hopf algebra $H$. The two
constructions are closely related, for one can construct, at least if $H$ is finite dimensional, an
equivalence  
\[ D(H)-\Mod \simeq Z_1(H-\Mod) \]
of braided tensor categories. Cf.\ \cite{kassel}. Since it is not true that all tensor categories
arise from Hopf algebras, the construction of the braided center $Z_1(\2C)$ can be considered a
generalization of the Hopf algebraic quantum double. (One might also find the definition of
$Z_1(\2C)$ more natural than that of $D(H)$.)

6. If $\2C$ is $k$-linear, spherical or a $*$-category (=unitary category), the same holds for
$Z_1(\2C)$. Other properties are much harder to show. In situations where a tensor category $\2C$ is
not of the form $H-\Mod$ for some Hopf algebra $H$, it can actually be quite difficult to construct  
objects of $Z_1(\2C)$ different from $\11_{Z_1(\2C)}$. There are situations where $\2C$ is quite
big, but $Z_1(\2C)$ is `trivial' (in the sense that $X\cong\11\ \forall X$): This happens if $\2C_0$
is a category and $\2C_1=Z_0(\2C_0)=\End(\2C_0)$, categorifying the simple fact that the center (in
the usual sense) $Z_1(M)$ of the monoid $M=Z_0(S)$ of endomaps of any set $S$ is trivial in the
sense of $Z_1(Z_0(S))=\{\id_S\}$. 

As we will see later, the situation is much better if $\2C$ is a fusion category.
\erem

In view of its construction, it is clear that there is a forgetful tensor functor 
\[ K: Z_1(\2C)\rarr \2C, \quad (X,e_X)\mapsto X. \]
In general, there is no natural functor, in particular not an inclusion, from $\2C$ to
$Z_1(\2C)$. There are two exceptions: If $\2C$ is a fusion category, $K$ has a 2-sided
(non-monoidal) adjoint $I:\2C\rarr Z_1(\2C)$. This will be discussed later.

The other exception is the case where $\2C$ comes with a braiding $c$. While the definition of
$Z_1(\2C)$ makes no reference to $c$, its existence has many consequences, cf.\ \cite{mue10}:

\begin{enumerate}
\item There are two tensor functors $F_1,F_2:\2C\rarr Z_1(\2C)$, given by
\[ F_1(X)=(X,e_X) \ \ \mbox{with}\ \ e_X(Y)=c_{X,Y}, \quad\quad
   F_2(X)=(X,\widetilde{e}_X) \ \ \mbox{with}\ \ \widetilde{e}_X(Y)=\widetilde{c}_{X,Y}.\]
Both functors are braided tensor functors from $\2C=(\2C,c)$ and
$\widetilde{\2C}=(\2C,\widetilde{c})$, respectively,  to $Z_1(\2C)$.
\item $F_1$ and $F_2$ are full and faithful, i.e.\ give embeddings 
$\2C\hookrightarrow Z_1(\2C),\ \widetilde{\2C}\hookrightarrow Z_1(\2C)$.
E.g.,
\bean \lefteqn{\Hom_{Z_1(\2C)}(F_1(X),F_1(Y)) = \Hom_{Z_1(\2C)}((X,c_{X,\bullet}),(Y,c_{Y,\bullet})) }\\
  && = \{ t\in\Hom_\2C(X,Y) \ | \  \id_X\otimes t \mcirc c_{X,Z}=c_{Y,Z}\mcirc t\otimes\id_X \quad
  \forall Z\in\2C \} =\Hom_\2C(X,Y), \eean
due to the naturality of $c$ w.r.t.\ both arguments.
\item The full subcategories $F_1(\2C)\subset Z_1(\2C)$ and $F_2(\widetilde{\2C})\subset Z_1(\2C)$
 commute with each other: By the definitions of $Z_1(\2C)$ and of $F_1,F_2$, we have
\[ c_{F_1(X),F_2(Y)}\circ c_{F_2(Y),F_1(X)}=c_{X,Y}\circ \widetilde{c}_{Y,X}=c_{X,Y}\circ
   c_{X,Y}^{-1}=\id_{Y\otimes X}. \]
An almost equally simple argument shows that $F_1(\2C), F_2(\widetilde{\2C})$ are each others
centralizers:  
\[ F_1(\2C)'=C_{Z_1(\2C)}(F_1(\2C))=F_2(\widetilde{\2C}), \quad\quad 
    F_2(\widetilde{\2C})'=C_{Z_1(\2C)}(F_2(\widetilde{\2C}))=F_1(\2C).\]
Consequentially, 
\[ F_1(\2C)''=F_1(\2C), \quad\quad F_2(\widetilde{\2C})''=F_2(\widetilde{\2C}).\]
This is nice since a priori we only know that $\2D\subset\2D''$ for a tensor subcategory
$\2D\subset\2C$. Again, for a fusion category $\2C$, one can prove much stronger results.
\item Since $F_1(\2C)$ and $F_2(\widetilde{\2C})$ are full subcategories of $Z_1(\2C)$, so is their
  intersection, and one finds
\[ F_1(\2C)\cap F_2(\widetilde{\2C})=F_1(Z_2(\2C))=F_2(Z_2(\2C)). \]
\item In view of item 4, $F_1$ and $F_2$ combine to a braided tensor functor
\[ H: \2C\times\widetilde{\2C}\rarr Z_1(\2C), \quad\quad (X,Y)\mapsto F_1(X)\otimes F_2(Y).\]
In view of 5., this functor will be neither full nor faithful in general: 
If $X\in Z_2(\2C)$ and $X\not\cong\11$ then $(X,\11)\not\cong(\11,X)$ but
$H((X,\11))=F_1(X)=F_2(X)=H((\11,X))$, thus $H$ does not reflect isomorphisms. However, we will see
that the linearized version $\2C\boxtimes\widetilde{\2C}\rarr Z_1(\2C)$ of $H$ actually is an 
equivalence when $\2C$ is modular.
\end{enumerate}

Another definition, due to \cite{bez}, involving the braided center will be useful later:

\bdefin \label{def-cs}
Let $\2C$ be a BTC, $\2D$ a tensor category and $F:\2C\rarr\2D$ a tensor functor. A central
structure on $F$ is a braided tensor functor $\widehat{F}:\2C\rarr Z_1(\2D)$ such that
$K\circ\widehat{F}=F$. Here $K:Z_1(\2D)\rarr\2D$ is the tensor functor that forgets the
half-braiding, thus $\widehat{F}$ is a lift of $F$ from $\2C$ to $Z_1(\2C)$.
\edefin

\brem From the point of view of the Baez-Dolan picture of `$k$-tuply monoidal $n$-categories',
cf.\ \cite{baez}, it is interesting to note the close analogy between central functors 
$\widehat{F}: \2C\rarr Z_1(\2D)$ (where $\2C$ is braided and $\2D$ just monoidal) and actions $F:
\2C\rarr Z_0(\2D)$ (with $\2D$ a category and $\2C$ monoidal). In both cases, the center  
$Z_0$ resp.\ $Z_1$ on the r.h.s.\ serves to create the piece of structure (monoidal,  braiding) that
is needed in order to talk about a  monoidal functor $F$ or braided functor $\widehat{F}$.
\erem

\subsection{Algebras and modules in BTCs and module categories of BTCs}\label{ss-algebras2}
\bdefin Let $\2C$ be a braided tensor category. An algebra $(A,m,\eta)$ in $\2C$ is called
commutative if $m\circ c_{A,A}=m$. A commutative separable algebra is called \'etale.
\edefin

In the case where $\2C$ is braided and the algebra $(A,m,\eta)$ in $\2C$ is commutative,
we would like ${}_A\2C$ to be a (braided) tensor category. In order for this to hold we need an
additional assumption: 

{\bf From now on, we require without further mention that $\2C$ has coequalizers.} (Cf.\ \cite{cwm}
for the definition.) 

Later on, all categories we consider will be at least abelian and therefore satisfy this assumption.

\bdefprop Let $\2C$ be a BTC and $(A,m,\eta)$ a commutative algebra, and let $(X,\mu),(X',\mu')\in{}_A\2C$.
Let $\alpha:X\otimes X'\rarr X''$ be a coequalizer of the pair of morphisms
\[ \mu_1=\mu\otimes\id_{X'}, \ \ \mu_2=(\id_X\otimes\mu')\circ(c_{A,X}\otimes\id_{X'}): \ 
  A\otimes X\otimes X'\rarr X\otimes X'.\]
By the universal property of $\alpha$, there is a unique $\mu'':A\otimes X''\rarr X''$ such that
$\mu''\circ(\id_A\otimes\alpha)=\alpha\circ\mu_1=\alpha\circ\mu_2$. Now $(X'',\mu'')\in{}_A\2C$.
With $(X,\mu)\otimes(X',\mu'):=(X'',\mu'')$, ${}_A\2C$ is a tensor category. (Commutativity of $A$ is
needed for the interchange law $(s\otimes t)\circ(s'\otimes t')=(s\circ s')\otimes(t\circ t')$ in
${}_A\2C$, to wit functoriality of $\otimes$ on morphisms.)
\edefprop

\blemma If $\2C$ is braided and $(A,m,\eta)$ a commutative algebra in $\2C$ then ${}_A\2C$ is a
tensor category (with the above tensor product). The free module functor $F_A:\2C\rarr{}_A\2C$ is a
tensor functor. 
\elemma

We would like to prove that ${}_A\2C$ is braided or symmetric. This requires additional
assumptions. If the algebra $(A,m,\eta)$ is \'etale, one can prove the following: 

\begin{itemize}
\item[(i)] If $\2C$ is symmetric then ${}_A\2C$ is symmetric and the functor $F_A$ is symmetric.
\item[(ii)] If $A\in Z_2(\2C)$ then ${}_A\2C$ is braided and the functor $F_A$ is braided. 
\item[(iii)] If $A\in Z_2(\2C)$ does not hold, the tensor category ${}_A\2C$ does not admit a
  braiding for which $F_A:\2C\rarr{}_A\2C$ is braided. The reason is that, every object of ${}_A\2C$
  being a quotient of $F_A(X)$ for some $X$, the only possible candidate for a braiding on
  ${}_A\2C$ making $F_A$ braided is the push-forward `$F_A(c)$' of the braiding $c$ of
  $\2C$. However, when $A\not\in Z_2(\2C)$, the would-be braiding $F_A(c)$ is natural only w.r.t.\ one of
  its arguments. Reformulating this positively, one obtains:
\item[(iv)] $F_A$ can always be considered as a braided tensor functor $\2C\rarr Z_1({}_A\2C)$. More
  precisely, the tensor functor 
  $F_A:\2C\rarr{}_A\2C$ admits a central structure $\widehat{F_A}:\2C\rarr Z_1({}_A\2C)$ in the
  sense of Definition \ref{def-cs}. Cf.\    \cite{dgno2,dmno}.
\item[(v)] An $A$-module $(X,\mu)$ is called {\it dyslectic} \cite{pareigis}, cf.\ also \cite{ko},
  or {\it local} when 
$\mu\circ c_{X,A}=\mu\circ\widetilde{c}_{X,A}$. The full subcategory ${}_A\2C^0\subset{}_A\2C$ of
  dyslectic modules is monoidal and in fact inherits a braiding from $\2C$. (Notice that
  ${}_A\2C^0={}_A\2C$ when $A\in Z_2(\2C)$, thus in particular when $\2C$ is symmetric.) The BTC
  ${}_A\2C^0$ will play an important r\^ole in the sequel.
\item[(vi)] In order to define a monoidal structure on ${}_A\2C$, where $A$ is an algebra in $\2C$,
  one actually does not need a braiding on all of $\2C$. Reviewing how the tensor structure was
  defined above, one realizes that   it suffices to be able to commute $A$ with all objects of
  $\2C$. More precisely, one should have a commutative algebra $((A,e_A),m,\eta)$ in $Z_1(\2C)$ ! In this
  situation, one has a  natural monoidal structure on ${}_A\2C$, cf.\ \cite{schauen}.  (When $\2C$
  is braided and $A\in\2C\stackrel{F_1}{\hookrightarrow} Z_1(\2C)$, this monoidal structure
  coincides with the one above since then $F_1(A)=(A,c_{A,\bullet})$.)
  Furthermore, for a commutative algebra $A\in Z_1(\2C)$, Schauenburg proved \cite{schauen} the
  remarkable braided equivalence   $Z_1({}_A\2C)\simeq {}_AZ_1(\2C)^0$.  
\end{itemize}

\subsection{Duality in braided categories. Braided fusion categories}\label{ss-db}
Let $\2C$ be a tensor category equipped with a left duality $X\mapsto({}^\vee X,e_X,d_X)$ and a
braiding $c$. Defining
\be e'_X=e_X\circ c_{X,{}^\vee X}: X\otimes{}^\vee X\rarr\11,\quad\quad\quad 
   d'_X=(c_{X,{}^\vee X})^{-1}\circ d_X: \11\rarr{}^\vee X\otimes X, \label{eq-ed'}\ee
a computation shows that $({}^\vee X,e'_X,d'_X)$ is a right dual for $X$. Thus in a braided category,
left and right duals of each object are isomorphic, and 
\be \alpha_X=(\id_{{}^{\vee\vee}X}\otimes e_X)\circ(d'_X\otimes\id_X): \ X\stackrel{\cong}{\rarr}{}^{\vee\vee}X
\label{eq-alpha}\ee
defines a natural isomorphism $\id\cong{}^{\vee\vee}-$ of functors. If $\2C$ is symmetric a computation
shows that $\alpha$ is a monoidal natural isomorphism, thus $\2C$ is pivotal. Another computation
shows that left and right traces coincide, thus $\2C$ is spherical.

All this breaks down if $\2C$ is braided but not symmetric, thus simply defining a right duality in terms
of a given left duality and the braiding does not give a satisfactory result. One solution is to
require in addition the existence of a ribbon structure:

\bdefin Let $\2C$ be a (strict) tensor category with braiding $c$ and left duality 
$X\mapsto({}^\vee X,e_X,d_X)$. A ribbon structure on $\2C$ is a natural isomorphism
$\Theta:\id_\2C\rarr\id_\2C$, i.e.\ a natural family of isomorphisms $\Theta_X: X\rarr X$,
satisfying
\bea \Theta_{X\otimes Y} &=& (\Theta_X\otimes\Theta_Y)\circ c_{Y,X}\circ c_{X,Y}\quad\forall X,Y \label{ribbon1}\\
 \Theta_{{}^\vee\!X} &=& {}^\vee(\Theta_X)\quad\forall X. \label{ribbon2}\eea
\edefin
(Notice that when the braiding $c$ is symmetric, (\ref{ribbon1}) makes the natural isomorphism $\Theta$ monoidal.)
Using the ribbon structure, we modify the formulas (\ref{eq-ed'}) as follows:
\[ e'_X=e_X\circ c_{X,{}^\vee X}\circ(\Theta_X\otimes\id_{{}^\vee X}): X\otimes{}^\vee X\rarr\11,
    \quad\quad\quad 
   d'_X=(\id_{{}^\vee X}\otimes\Theta_X)\circ(c_{X,{}^\vee X})^{-1}\circ d_X: \11\rarr{}^\vee X\otimes X.\]
Now one finds that $\{ \alpha_X\}$, defined as in (\ref{eq-alpha}), but using the modified
definitions of $d'_X, e'_X$, is a monoidal natural isomorphism. Thus $\2C$ is pivotal, and again in
fact spherical. Cf.\ \cite{kassel}.

\brem Occasionally, it is preferable to reverse the above logic. Namely, if $\2C$ is a spherical
category and $c$ a braiding (with no compatibility assumed) then defining 
\[ \Theta_X= (tr_X\otimes\id_X)(c_{X,X}), \]
one finds that $\{ \Theta_X\}$ satisfies (\ref{ribbon1},\ref{ribbon2}), thus is a ribbon structure
compatible with the braiding $c$. Furthermore, the natural isomorphism $\alpha:\id\rarr{}^{\vee\vee}-$ given
as part of the spherical structure coincides with the one defined in terms of the left duality and
$\Theta$ as in (\ref{eq-alpha}). Therefore, for a braided category $(\2C,c)$ with left duality,
giving a pivotal (in fact spherical) structure $\alpha$ is equivalent to giving a ribbon structure
$\Theta$ compatible with $c$. Cf.\  \cite{yetter1}.
\erem

We now briefly return to the subject of algebras in braided categories and their module
categories. By a braided fusion category we simply mean a fusion category equipped with a braiding,
and similarly for braided spherical categories.  Now one has:

\bprop \label{prop-CA}\cite{ko}
If $\2C$ is a braided fusion (resp.\ spherical) category and $A\in\2C$ an \'etale algebra
then ${}_A\2C$ is a fusion (resp.\ spherical) category (not necessarily braided) and
\[ \mathrm{FP}\dim {}_A\2C=\frac{\mathrm{FP}\dim\2C}{d(A)}. \]
If $\2C$ is spherical, both instances of FP$\dim$ in this identity can be replaced by $\dim$ as
defined in terms of the spherical structures.
\eprop

The dimension of the braided fusion category ${}_A\2C^0$ does not just depend on (FP)$\dim\2C$ and
$d(A)$, but on `how much' of the object $A$ lies in $Z_2(\2C)$. (E.g., it is evident that
${}_A\2C^0={}_A\2C$ when $A\in Z_2(\2C)$, giving the same dimension for both categories. But as we
will see, it is also possible that ${}_A\2C^0$ is trivial, i.e.\ has dimension $1$.) However, when $Z_2(\2C)$ is trivial,
i.e.\ $\2C$ is non-degenerate (resp.\ modular), one again has a simple formula, cf.\ (\ref{eq-dim0}) below.  

As discussed earlier, the tensor functor $F_A:\2C\rarr{}_A\2C$, while not braided in general, always
admits a central structure. In the setting of fusion categories, central functors and module
categories are closely related:

\btheor \cite[Lemmas 3.5, 3.9]{dmno} (i) If $\2C$ is a braided fusion category, $\2D$ is a fusion category and $F:
\2C\rarr\2D$ a central functor then there exists a connected \'etale algebra $A\in\2C$
such that the category ${}_A\2C$  is monoidally equivalent to the image of $F$, i.e.\ the 
smallest fusion subcategory of $\2D$ containing $F(\2C)$. (The object $A$ is determined by
$\Hom_\2C(X,A)\cong\Hom_\2D(F(X),\11)$ and exists since $F$ has a right adjoint.)

(ii) If $\2C$ is a braided fusion category and $A\in\2C$ a connected \'etale algebra then the
connected \'etale algebra $A'\in\2C$ obtained by (i) from the central functor $F_A:\2C\rarr{}_A\2C$ is isomorphic to
$A$.
\etheor


\section{Modular categories}
\subsection{Basics}
The rest of this paper will be concerned with braided fusion categories over an algebraically closed
field $k$, most often $\7C$. Recall that every fusion category has a minimal fusion subcategory
consisting only of the multiples of the unit $\11$. This subcategory is equivalent to Vect$_k$. A
fusion category $\2C$ is called trivial when it is itself equivalent to 
Vect$_k$, which is the same as saying that every simple $X\in\2C$ is isomorphic to $\11$.

\bdefin A braided fusion category $\2C$ is called
\begin{itemize}
\item pre-modular if it is spherical,
\item non-degenerate if $Z_2(\2C)$ is trivial,
\item modular if it is pre-modular and non-degenerate. Such a $\2C$ is called just `modular
  category' or MTC.
\end{itemize}
\edefin

\brem 1. Non-degenerate braided fusion categories are related to symmetric fusion categories like
von Neumann factors, i.e.\ von Neumann algebras $M$ with trivial center $Z(M)$, to commutative von
Neumann algebras, where $Z(M)=M$. Since these two extremal types of von Neumann algebras play
distinguished r\^oles in the general theory, it should not come as a surprise that the analogue also
holds in the setting of braided fusion categories. 

2. By an important theorem of Doplicher and Roberts \cite{DR} and independently Deligne \cite{del},
cf.\ Section \ref{s-modularization}, symmetric
fusion categories are closely related to finite groups (and supergroups). Thus classifying symmetric
fusion categories is essentially equivalent to classifying finite groups, a rather difficult task that
has been achieved only partially. Given the importance of modular categories in the contexts
of quantum group theory, conformal field theory, low dimensional topology, in particular topological
quantum field theories, one may argue that the study and classification of modular categories is 
as natural and urgent as that of finite groups.
\erem

Let $\2C$ be a pre-modular category. For $X,Y\in\2C$, define
\[ S(X,Y)=\tr_{X\otimes Y}(c_{Y,X}\circ c_{X,Y})\ \in k.\]
By the properties of the trace, $S(X,Y)$ depends only on the isomorphism classes $[X],[Y]$. Thus if 
$I(\2C)$ is the set of isomorphism classes of simple objects of $\2C$ and we choose representers
$X_i$, we can define an $|I(\2C)|\times |I(\2C)|$-matrix $S$ by $S_{i,j}=S(X_i,X_j)$. $S$ is
symmetric, and it is easy to see that non-triviality of $Z_2(\2C)$ implies singularity of $S$: If
$X_i\in Z_2(\2C)$ then $S_{i,j}=\tr_{X_i\otimes X_j}(\id)=d(X_i)d(X_j)$ for all $j$, and thus the
$i$-th row (and column) are proportional to the $0$-th row (column) (where $X_0\cong\11$).

More interestingly, for $\2C$ pre-modular one can show, cf.\ \cite{turaev1,turaev,khr}:
\begin{itemize}
\item For simple $X,Y$, one has $S(X,Y)=d(X)d(Y)$  if \underline{and only if} $X$ and $Y$ commute.
\item Let $\2K\subset\2C$ be a fusion subcategory. Then
\[ \sum_{i\in I(\2K)} d(X_i) S_{i,j}=\left\{ \begin{array}{cc} d(X_j)\dim\2K & \mbox{if}\ X_j\in \2K'\\
   0 & \mbox{otherwise} \end{array}\right. \]
\item If $Z_2(\2C)$ is trivial then $S^2=\dim\2C\,C$, where $C$ is the `charge-conjugation' matrix:
  $C_{i,j}=\delta_{i,\ol{\jmath}}$. (Note that $C^2=\11$.) Thus if we assume $\dim\2C\ne 0$ (which
  is automatic over $\7C$ \cite{eno}), then $S$ is invertible if and only if $Z_2(\2C)$ is trivial,
  i.e.\ $\2C$ is modular. 
\item If $\2C$ is modular then the `Gauss sums' of $\2C$, defined by 
\begin{equation}\label{eq-gauss} \Omega^\pm(\2C)=\sum_{i\in I(\2C)} \Theta(X_i)^{\pm 1} d(X_i)^2 \ee 
satisfy $\Omega^+(\2C) \Omega^-(\2C)=\dim\2C$ and the diagonal $|I(\2C)|\times |I(\2C)|$-matrix $T$
defined by   $T_{i,j}=\delta_{i,j}\Theta(X_i)$ satisfies $TSTST=\Omega^+(\2C) S$.
(The modular category $\2C$ is called `anomaly-free' when $\Omega^+=\Omega^-$.)
\item Therefore, if $\2C$ is modular, thus $S$ invertible, we have $S^2=\alpha C$ and $(ST)^3=\beta C$
  with $\alpha\beta\ne 0$. This means that $S$ and $T$ define a projective representation of the
  modular group $SL(2,\7Z)$. (Recall that the latter is generated by the elements 
\[ s=\left( \begin{array}{cc} 0 & 1 \\ -1 & 0 \end{array}\right), \quad\quad t=\left( \begin{array}{cc} 1 & 1 \\ 0 & 1 \end{array}\right), \]
which satisfy $s^2=-\11, (st)^3=\11$.) The existence of this representation is the rationale behind the
terminology `modular categories'.
\end{itemize}

Let $X$ be a simple object and $Y$ invertible. Then $X\otimes Y$ is simple, thus
$c_{Y,X}\circ c_{X,Y}\in\Aut(X\otimes Y)$ is a scalar (=element of the ground field). Multiplying
this scalar by $d(X)$ gives $S(X,Y)$, but this is not the point. The point is that the map
$f:I\times I_1\rarr k$, where $I_1\subset I$ is the subgroup of invertible isomorphism classes,
obtained by the above consideration is a homomorphism  w.r.t.\ the second argument. One can show
that $f(Z,L)=f(X,L)f(Y,L)$ whenever $N_{X,Y}^Z>0$, cf.\ \cite[Section 4]{khr}. This implies
$f(X,Z)=f(Y,Z)$ whenever $\del X=\del Y$,   where $\del: I\rarr G(I)$ is  the universal group
grading discussed in Section \ref{ss-fusion}. Thus $f$ descends to a bihomomorphism 
$G(I)\times I_1\rarr k$. Remarkably one has:

\btheor \cite{GN1}\label{theor-GN}
When $\2C$ is modular, the above map is a non-degenerate pairing, establishing a 
canonical isomorphism   $I_1\rarr\widehat{G(I)}$.
\etheor

\brem Recall that for a finite group $H$, the grading group of $\Rep\,H$ is given by
$G(\widehat{H})\cong\widehat{Z(H)}$. On 
the other hand, $I_1(\Rep\,H)\cong \widehat{H_{ab}}$, where $H_{ab}=H/[H,H]$ is the abelianization. 
The abelian groups $Z(H)$ and $H_{ab}$ have little to do with each other. In view of this, Theorem
\ref{theor-GN} is one of many manifestations of the observation of \cite{khr} that a modular
category is ``a self-dual object that is more symmetric than a group''.
\erem

In the next two sections, we will encounter several constructions that give rise to modular
categories. However, it seems instructive to give an example already at this point. 

\bexam\label{ex}
Let $A$ be a finite abelian group and $\widehat{A}$ its character group. Let $\2C_0(A)$ be the
strict tensor category defined by  
\begin{itemize}
\item[(i)] $\obj\,\2C(A)=A\times\widehat{A}$.
\item[(ii)] $\Hom((g,\phi),(g',\phi'))=\7C$ if $(g,\phi)=(g',\phi')$ and $=\{0\}$ otherwise.
\item[(iii)] $(g,\phi)\otimes(g',\phi')=(gg',\phi\phi')$. Composition and tensor product of
  morphisms are defined as multiplication of complex numbers.
\item[(iv)] The braiding is given by $c_{(g,\phi),(g',\phi')}=\phi(g')\id$. (This makes sense since $(g,\phi)\otimes(g',\phi')=(g',\phi')\otimes(g,\phi)$.)
\end{itemize}
Now $\2C(A)$ is the closure of $\2C_0(A)$ w.r.t.\ direct sums. One finds $\Theta((g,\phi))=\phi(g)$ and
$S((g,\phi),(g',\phi'))=\phi(g')\phi'(g)$, from which it follows easily that $(e,\phi_0)$, where
$\phi_0\equiv 1$, is the only central object, thus $\2C(A)$ is modular.

The above construction of $\2C(A)$ admits generalization to non-abelian finite groups, but that is
better done using Hopf algebra language. This leads to the quantum double $D(G)$ of a finite group
or a finite dimensional Hopf algebra $D(H)$, cf.\ the first paragraph of Section \ref{s-z1}. 
\eexam

\bitem For any finite abelian group, one easily proves that $\Omega^\pm(\2C(A))=N$. Thus the
categories $\2C(A)$ are anomaly-free. If $N$ is odd then the full subcategory
$\2D_N\subset\2C(\7Z/N\7Z)$ with objects  
$\{ (k,k) \ | \ k=0,\ldots,N-1\}$ is itself modular, 
cf.\ Remark \ref{rem-fact}.3, and one has
\[ \Omega^\pm(\2D_N)=\sum_{k=0}^{N-1} e^{\pm\,\frac{2\pi i}{N}k^2}, \]
which is a classical Gauss sum. (This motivates calling the quantities $\Omega^\pm(\2C)$ `Gauss
sums' in general.) By the classical computation of these Gauss sums, we find that
$\Omega^+(\2D_N)$ equals $\sqrt{N}$ when  $N\equiv 1(\mbox{mod}\ 4)$ and $i\sqrt{N}$ when $N\equiv
3(\mbox{mod}\ 4)$. In view of 
$\Omega^-(\2C)=\ol{\Omega^+(\2C)}$, the categories $\2D_N$ with $N\equiv\ 3(\mbox{mod}\ 4)$
therefore provide our first examples of modular categories that are not anomaly-free.
\eitem

\bitem
Modular categories first appeared explicitly in \cite{turaev1}, but some aspects of that paper
were anticipated by two years in \cite{khr}, which also was inspired by \cite{ms1} but in addition
drew upon the well established operator algebraic approach \cite{haag} to axiomatic quantum field
theory and in particular on \cite{frs}. In particular, \cite{khr} contains the first proof of the
equivalence between invertibility of the $S$-matrix, in terms of which modularity was defined in
\cite{turaev1,turaev}, and triviality of $Z_2(\2C)$. (For a more general recent proof cf.\ \cite{bebl2}.)
\eitem


\subsection{Digression: Modular categories in topology and mathematical physics}
\bitem
While some {\it inspiration} for \cite{turaev1} came from conformal field theory and in particular
\cite{ms1}, cf.\ below, the main {\it motivation} came from 
low-dimensional topology. In 1988/9, Witten \cite{wit} had proposed an interpretation of the new
`quantum invariants' of knots and 3-manifolds (in particular Jones' polynomial knot invariant) in
terms of `topological quantum field theories' 
(TQFTs), defined via a non-rigorous (not-yet-rigorous?) path-integral formalism. From this work,
Atiyah \cite{ati} immediately abstracted the mathematical axioms that should be satisfied by a 
TQFT, and Reshetikhin and Turaev \cite{rt2} soon used the representation theory of quantum
groups to rigorously construct a TQFT that should essentially be that studied by Witten. The aim of
\cite{turaev1} was then to isolate the mathematical structure that is behind the construction in
\cite{rt2} as to enable generalizations, and indeed one has a $2+1$-dimensional TQFT $F_\2C$ for
each modular category $\2C$, wherever $\2C$ may come from. (Cf.\ \cite{turaev} for a full exposition
of the early work on the subject and \cite{BK} for a somewhat more recent introduction.)

Since there is no space for going into this subject to any depth, we limit ourselves to mentioning
that a TQFT in $2+1$ dimensions gives rise to projective representations of the mapping class
groups of all closed two-manifolds. The mapping class group of the 2-torus is the modular group
$SL(2,\7Z)$, and its projective representation produced by the TQFT $F_\2C$ associated with the
modular category $\2C$ is just the one encountered above in terms of $S$ and $T$.
\eitem

\bitem
Before we turn to a brief discussion of the r\^ole of modular categories in conformal field theory,
we mention another manifestation of them in mathematical physics. It has turned out that
infinite quantum systems (both field theories and spin systems) in two spatial dimensions can have
`topological excitations', whose mathematical analysis leads to braided tensor categories which
often turn out to be modular. An important example is Kitaev's `toric code' \cite{kitaev}, a quantum
spin system in two dimensions which gives rise to the (rather simple) modular category
$D(\7Z/2\7Z)-\Mod$ (the $\2C(\7Z/2\7Z)$ of Example \ref{ex}). (Cf.\ also \cite{naaijkens}.) There is
a generalization of the toric code to finite groups other than $\7Z/2\7Z$, but not everything has
been worked out yet. The toric code models play a prominent r\^ole in the subject of topological
quantum computing, reviewed e.g.\ in \cite{freedman2, freedman1, wang}. 
\eitem

\bitem
While Kitaev's model lives in $2+1$ dimensions and has a mass gap, braided and modular categories
also arise in conformally invariant (thus massless) quantum field theories in $1+1$ or $2+0$
dimensions, a subject that has been researched very extensively. In particular, there have been two
rigorous and model-independent proofs of the statement that suitable chiral conformal field theories
have modular representation theories. In the operator algebraic approach to CFTs, the basic
definitions are quite easy to state, and we briefly do so. 

The group $\2P=PSU(1,1)$ acts on $S^1=\{ z\in\7Z\ | \ |z|=1\}$ by
$\left(\begin{array}{cc} a & b\\ c & d\end{array}\right)z=\frac{az+b}{cz+d}$.
Let $\2I$ be the set of connected open subsets $I\subset S^1$ such that $\emptyset\ne I\ne
S^1$. (Thus $\2I$ consists of the non-trivial connected open intervals in $S^1$.) If $I\in\2I$ and
$g\in\2P$ then $gI\in\2I$.

\bdefin A chiral CFT $\2A$ is a quadruple $(H_0,A(\cdot),U,\Omega)$, where $H_0$ is a separable Hilbert
space (thus essentially independent of $\2A$), $\Omega\in H$ a unit vector, $U$ is a strongly
continuous positive-energy representation of the group $PSU(1,1)$ on $H$, and for each $I\in\2I$,
$A(I)\subset B(H_0)$ is a von Neumann algebra. These data must satisfy the following axioms:
\begin{itemize}
\item[(i)] `Isotony': $I_1\subset I_2\impl A(I_1)\subset A(I_2)$
\item[(i)] `Locality': $I_1\cap I_2=\emptyset\impl [A(I_1),A(I_2)]=\{0\}$.
\item[(iii)] `Irreducibility': $\cap_I A(I)'=\7C\11$.
\item[(iv)] `M\"obius Covariance': If $I\in\2I, g\in\2P$ then $U(g)A(I)U(g)^*=A(gI)$.
\item[(v)] $U(g)\Omega=\Omega\ \forall g\in\2P$.
\end{itemize}
\edefin
Already from these few axioms one can prove a number of important results, among which:

\begin{itemize}
\item[(a)] `Factoriality': Each $A(I)$ is a factor, i.e.\ has trivial center. 
\item[(b)] `Weak additivity': If $I_1,I_2\in\2I$ such that $I_1\cup I_2\in\2I$ then 
$A(I_1)\vee A(I_2)=A(I_1\cup I_2)$ . 
\item[(c)] `Haag duality': For each $I\in\2I$ one has $A(I)'=A(I')$, where $I'$ is the interior of $S^1-I$.
\end{itemize}
One has a notion of representation for chiral CFTs:

\bdefin Let $\2A=(H,A(\cdot),U,\Omega)$ be a chiral CFT. A representation of $\2A$ is a pair
$(H,\{\pi_I\}_{I\in\2I})$, where $H$ is a (separable) Hilbert space and
$\{\pi_I: A(I)\rarr B(H)\}_{I\in\2I}$ is a family of $*$-representations satisfying 
$\pi_{I_2}\restr A(I_1)=\pi_{I_1}$ whenever $I_1\subset I_2$.
A morphism of representations $(H,\{ \pi_I\}), (H',\{ \pi'_I\})$ is a bounded operator $V:H\rarr H'$
such that $V\pi_I(\cdot)=\pi'_I(\cdot)V$ for each $I\in\2I$. This defines a $*$-category
$\Rep\,\2A$. 

If $(H,\{\pi_I\})$ is a representation, we define $\dim\pi=[\pi_{I'}(A(I')'):\pi_I(A(I))]\in[1,\infty]$, 
where the brackets denote the Jones index, and the r.h.s.\ is independent of $I\in\2I$. $\Rep_f\2A$
denotes the full subcategory of $\Rep\,\2A$ of representations with $\dim\pi<\infty$. 
\edefin

It is clear that $\Rep\,\2A$ and $\Rep_f\2A$ are $*$-categories. But one prove much more: Both admit
canonical braided monoidal structures, and $\Rep_f\2A$ is semisimple with (two-sided) duals, but not
necessarily finite. Cf.\ \cite[II]{frs}. In order to prove finiteness or modularity of $\Rep_F\2A$
one needs to assume more: 

\bdefin \label{def-klm}
A chiral CFT $\2A=(H,A(\cdot),U,\Omega)$ is called completely rational if the following holds:
\begin{itemize}
\item[(i)] `Split property': If $I_1,I_2\in\2I$ satisfy $\ol{I_1}\cap\ol{I_2}=\emptyset$ then the
  natural map $A(I_1)\otimes_{alg} A(I_2)\rarr A(I_1)\vee A(I_1)$ induces an isomorphism of von
  Neumann algebras. 
\item[(ii)] `Strong additivity': If $I_1,I_2\in\2I$ satisfy $\#(\ol{I_1}\cap\ol{I_2})=1$ (i.e.\ the
  intervals are disjoint but adjacent) then $A(I_1)\vee A(I_2)=A(I)$, where $I$ is the interior of
  the closure of $I_1\cup I_2$. 
\item[(iii)] `Finiteness': Let $I_1,I_2\in\2I$ such that $\ol{I_1}\cap\ol{I_2}=\emptyset$. Then
  $(I_1\cup I_2)'=I_3\cup I_4$ for certain $I_3,I_4\in\2I$. The quantity 
  $\mu(A)=[(A(I_3)\vee A(I_4)')':A(I_1)\vee A(I_2)]$ 
  (which a priori is in $[1,\infty]$ and can be shown to be independent of $I_1,I_2$) is finite.
\end{itemize}
\edefin

(Axiom (ii) is not really essential.) The following was proven in \cite{klm}:

\btheor \label{theor-klm} If $\2A=(H,A(\cdot),U,\Omega)$ is a completely rational chiral CFT then
$\Rep_f\2A$ is  modular with $\dim\Rep_f\2A=\mu(A)$.
\etheor

Many examples of completely rational chiral CFTs are known, and Theorem \ref{theor-klm} is just the
beginning of a rapidly growing theory. (The fact that the dimension of the representation category
can (in principle) be obtained by looking at just a few local algebras is extremely useful.)
Cf.\ e.g.\ the references in \cite{mue-prag}.

Around 2005, Huang proved a similar result in the framework of vertex operator algebras,
cf.\ \cite{huang} and references therein, assuming the property of `$C_2$-cofiniteness', which is
similar to (iii) in Definition \ref{def-klm}.

The most important chiral conformal field theories are
related to (a) the projective `positive energy' representations of the diffeomorphism group
Diff$^+(S^1)$ of the circle and/or its (centrally extended) Lie algebra, the Virasoro algebra, and
(b) the positive energy representations of the loop groups $C^\infty(S^1,G)$, where $G$ is a compact
Lie group. These representation categories can be studied without reference to conformal field
theory, cf.\ e.g.\ \cite{PS,wakimoto}, but defining the tensor structure does require conformal
field theory (or at least techniques of the latter), like those of \cite{frs} or \cite{wass}.
\eitem

\bitem \label{renorm} For any algebraically closed field $k$, the group $H^2(SL(2,\7Z),k^*)$ is
trivial. This 
  implies that by rescaling the matrices $S,T$ one can obtain a true representation of the modular
  group, and there are exactly six ways of doing this. There are two situations where there is a
  preferred choice: When $\2C$ is anomaly-free, as when $\2C\simeq Z_1(\2D)$, one has
  $\Omega^+=\Omega^-=\pm\dim\2C$. Then the renormalization $S'=S/\Omega^+, T'=T$
  gives a true representation of $SL(2,\7Z)$, which may be   considered more canonical than the
  others. 

On the other hand, when the modular category arises as the representation category of a conformal
field theory $\2A$ (in the setting of operator algebras or vertex algebras), 
for each (equivalence classes) of
  simple objects of $\Rep\,\2A$, there is an analytic function $f_i:\7H\rarr\7C$, the `character' of
  that representation. Collecting these $n=|I(\Rep\,\2A)|$ functions in a vector-valued function
  $f:\7H\rarr V=\7C^n$, one finds that $f$ satisfies the following definition:

\bdefin Let $V$ be a finite dimensional complex vector space and $\pi: SL(2,\7Z)\rarr\End(V)$ a
representation. A vector valued modular form of type $\pi$ is a holomorphic map $\rho:\7H\rarr V$ satisfying

\be\label{eq-trf} f(g^{-1}z)=\pi(g)f(z)\quad\forall z\in\7H,\ g\in SL(2,\7Z). \ee
Here $\7H=\{ z\in\7C\ | \ \mathrm{Im}\,z>0\}$ and $SL(2,\7Z)$ acts on $\7H$ by
$\left(\begin{array}{cc}  a&b\\ c&d\end{array}\right)z=\frac{az+b}{cz+d}$. 
\edefin

Here, $\pi$ is a uniquely determined (by the CQFT) true representation of $SL(2,\7C)$ obtained by
a particular renormalization the matrices $S,T$ associated with the modular category $\Rep\,\2A$.

In the CQFTs associated with the representation theories of the Virasoro and the Kac-Moody algebras, 
the above vector valued modular forms can be studied very explicitly, cf.\ e.g.\ \cite{wakimoto}.
\eitem

\bitem In all rational chiral CFTs that have been studied explicitly, it turned out that there is an
$N\in\7N$ such that all conformal characters $\chi_i$ are modular functions for the congruence
subgroup $\Gamma(N)=\{ M\in SL(2,\7Z) \ | \ M\equiv\11\ (\mod\ N)\}$. This means that
$\chi_i(gz)=\chi_i(z)$ for all $g\in \Gamma(N)$ and all $i$. In view of (\ref{eq-trf}) this amounts
to $\Gamma(N)\subset\mathrm{ker}(\pi)$. This led to the folk conjecture that this `congruence
subgroup property' holds in all rational chiral CFTs. Such a general result was indeed obtained in
\cite{bantay} (where unfortunately no rigorous formalism of CFTs was used). The `conductor' $N$ is
closely related to the order of the diagonal matrix $T$ (all elements of which are roots of unity).
It is natural to ask whether a similar result can be proven for all modular categories irrespective
of whether they arise from a chiral CFT. (The latter question is an important open problem.)
The answer is yes, but before we state it, we revert back to `pure' mathematics.
\eitem


\subsection{Modular categories: Structure theory and module categories}
The work \cite{bantay} inspired \cite{SZ} where a congruence subgroup theorem was proven for the
modular categories of the form $D(H)-\Mod$ and finally the very recent \cite{ng3} with a result
valid for all modular categories:

\btheor \cite{ng3} (a) If $\2D$ is a spherical fusion category and $\2C=Z_1(\2D)$ then the kernel of
the canonical (cf.\ \ref{renorm}) true representation of $SL(2,\7Z)$ contains $\Gamma(N)$, where $N$
is the `Frobenius-Schur exponent' of $\2C$.

(b) If $\2C$ is an arbitrary modular category, then $\Gamma(N)$ is contained in the kernel of the
canonical projective representation of $SL(2,\7Z)$ generated by $S,T$.
\etheor

The following results from \cite{mue11}, generalized to non-degenerate braided fusion categories in
\cite{dgno}, are the first steps towards a structure theory and perhaps classification of modular
categories. Part (iv) shows that, in a sense, modular categories are better behaved than finite groups!

\btheor \label{theor-factoriz}
Let $\2C$ be a modular category and $\2D_1\subset\2C$ a fusion subcategory. Let
$\2D_2=C_\2C(\2D_1)$. Then
\begin{itemize}
\item[(i)] $\dim\2D_1\cdot\dim\2D_2=\dim\2C$.
\item[(ii)] $C_\2C(\2D_2)=\ol{\2D_1}$, where the r.h.s.\ denotes the closure of $\2D_1$ under isomorphisms
  (i.e.\ the smallest replete fusion subcategory of $\2C$ containing $\2D_1$). Thus if $\2D\subset\2C$
  is a replete fusion subcategory then $\2D''=\2D$. 
\item[(iii)] If $\2D\subset\2C$ is a full tensor subcategory then $Z_2(\2D)=Z_2(\2D')$. In particular, if
  $\2D$ is modular then so is $\2D'$.
\item[(iv)] If $\2D\subset\2C$ is modular then $\2C\simeq\2D\boxtimes\2D'$ as BTC.
\end{itemize}
\etheor

\bdefin A modular category $\2C$ is called prime if every modular fusion subcategory is either
trivial or equivalent to $\2C$.
\edefin

\bcoro Every modular category is equivalent to a finite direct product of prime modular categories.
\ecoro

\brem \label{rem-fact} 1. Statement (ii) is similar to the `double commutant theorem' in the theory
of von Neumann algebras: The commutant $M''$ of the 
commutant of a von Neumann algebra equals $M$. Statement (iv), according to which modular
subcategories always are direct factors, is analogous to the fact that an inclusion $N\subset M$ of
type I factors gives rise to an isomorphism $M\cong N\otimes N'$.

2. If $G$ is a finite non-abelian simple group then the modular category $D(G)-\Mod$ is prime. In
fact it contains only one fusion subcategory, namely $\Rep\,G$. It follows that the classification of
prime modular categories contains the classification of finite simple (non-abelian) groups! For
general finite groups, $D(G)-\Mod$ can have many fusion subcategories. Cf.\ \cite{nnw} for a
classification. 

3. In general, the prime factorization of a modular category is not unique. An example is already
provided by $\2C=D(A)-\Mod$, where $A=\7Z/p\7Z$ is cyclic of prime order $p\ne 2$. In this case, the
replete prime modular subcategories of $\2C$ are labeled by the isomorphisms
$\phi:A\rarr\widehat{A}$, and $(\2D_\phi)'=\2D_{\ol{\phi}}$ where
$\ol{\phi}(\cdot)=\phi(\cdot)^{-1}$. Cf.\ \cite{mue11}. In this example, all objects of $\2C$ are
invertible, which is crucial for the non-uniqueness: By \cite[Proposition 2.2]{dmno}, the prime
factorization of $\2C$ is unique if $\2C$ has no invertible objects apart from $\11$. More
generally, prime factors having no invertibles other than $\11$ appear identically in every prime
factorization of $\2C$. Thus the non-uniqueness results from the possibility of `moving invertible
objects from one direct factor to another'. The details have not yet been clarified, but it is clear
that homomorphisms $\phi: G(\2C)\rarr I_1(\2C)\cong \widehat{G(\2C)} $ play a r\^ole. (Recall from
\ref{theor-GN} that the abelian grading group  $G(\2C)$  and the group $I_1(\2C)$ of (isomorphism classes
of) invertible objects of $\2C$ are canonically dual to each other.)
\erem

If $\2C$ is modular and $\2D\subset\2C$ a fusion subcategory then 
$Z_2(\2D)=\2D\cap\2D'\subset\2C\cap\2D'=C_\2C(\2D)$, implying  
\[ \dim\2C=\dim\2D\cdot\dim C_\2C(\2D)\ge\dim\2D\cdot\dim Z_2(\2D). \]
This motivated the conjecture \cite{mue11}:

\bconj \label{conj-embed} If $\2D$ is a pre-modular category, there is a modular category $\2C$
containing $\2D$ as a fusion subcategory such that $\dim\2C=\dim\2D\cdot Z_2(\2D)$.
\econj
(There are indications that the conjecture in this generality may be false, but see Theorem
\ref{theor-4} below.)

Above, we pointed out the analogy between non-degenerate braided fusion categories and von Neumann
factors. This analogy goes a bit further: Since factors are simple (as algebras), their
homomorphisms have trivial kernels and therefore are embeddings. Analogously, one has:

\bprop \label{prop-faithful}\cite[Corollary 3.26]{dmno} Any braided tensor functor $F: \2C \to \2D$
between braided fusion categories with $\2C$ (almost) non-degenerate is fully faithful.
\eprop

\noindent{\it Sketch of Proof.} Replacing $\2D$ by its subcategory $F(\2C)$, we may assume that
$F$ is surjective. By results mentioned earlier, there is an \'etale algebra $A\in\2C$ and an
equivalence $H:\2D\rarr{}_A\2C$ of fusion categories such that $H\circ F\cong F_A$. That $F_A$ is
braided implies $A\in Z_2(\2C)$, and $\2C$ being non-degenerate, we have $A=\11$, thus $F\cong\id$.
\qed

In analogy to Proposition \ref{prop-CA}, one has the following result concerning dyslectic module
categories for \'etale algebras in non-degenerate braided fusion categories:

\btheor \cite{ko}
Let $\2C$ be a non-degenerate braided fusion category and $A\in\2C$ a connected
\'etale algebra. Then the dyslectic module category ${}_A\2C^0$ is a non-degenerate braided
fusion category and 
\be \mathrm{FP}\dim {}_A\2C^0=\frac{\mathrm{FP}\dim\2C}{d(A)^2}. \label{eq-dim0}\ee
If $\2C$ is spherical, so is ${}_A\2C^0$ and FP$\dim$ can be replaced by $\dim$.
\etheor

In the next two sections we discuss two ways of constructing modular categories (or, more generally,
non-degenerate braided fusion categories). The first, modularization, starts from a pre-modular
category and the second from a mere (spherical) fusion category.


\section{Modularization of pre-modular categories. Generalizations}\label{s-modularization}
Let $\2C$ be a pre-modular category, to wit a braided spherical fusion category. 
Its symmetric center $Z_2(\2C)$ then is a symmetric spherical fusion category.
Since non-modularity of $\2C$ is equivalent to $Z_2(\2C)$ being non-trivial, it is natural to try to 
`quotient out' the full subcategory $Z_2(\2C)$ in order to obtain a modular category
`$\2C/Z_2(\2C)$'. Formalizing this idea one arrives at the following:

\bdefin \cite{brug} A modularization of a pre-modular category $\2C$ is a functor $F:\2C\rarr\2D$ of
braided fusion categories, where $\2D$ is modular and $F$ is surjective (or `dominant') in the sense
that every object of $\2D$ is a subobject of one of the form $F(X)$ with $X\in\2C$. 
\edefin

The fact that $F$ is supposed to be braided and surjective implies that it must trivialize
$Z_2(\2C)$, i.e.\ $F(X)$ must be a multiple of the identity whenever $X\in Z_2(\2C)$. 
The following was shown in \cite{brug,mue06}:

\btheor A pre-modular category admits a  modularization if and only if the symmetric category
$Z_2(\2C)$ is purely even, i.e.\ all objects have trivial twist, i.e.\ $\theta_X=\id_X\ \forall X$. 
\etheor

The proof relies on the following deep result:

\btheor \label{theor-tann}
 \cite{DR,del} Let $\2C$ be a spherical symmetric fusion category with trivial twists. Then there is
a finite group $G$, unique up to isomorphism, such that $\2C\simeq \Rep\,G$ as symmetric tensor
category. This equivalence is unique up to natural monoidal isomorphism.
\etheor
(Both \cite{DR,del} prove much more general results. In \cite{DR}, $\2C$ is  required to be unitary.)
Now, if $G$ is a finite group, the vector space $A=$Fun$(G,\7C)$ underlying the left regular 
representation also is a commutative algebra, and one finds that $A$ is a connected \'etale algebra
in $\Rep\,G$. Furthermore, the module category ${}_A(\Rep\,G)$ is trivial. We call $A$ the regular
algebra, and we also do this for the corresponding object in $\2C$ when Theorem \ref{theor-tann} is
invoked. Combining this observation with Theorem \ref{theor-tann}, it is clear how to obtain a
modularization of a pre-modular category $\2C$ with trivial twists: Take $\2D={}_A\2C$, where $A$ is
the regular algebra of the symmetric even category $Z_2(\2C)\subset\2C$, and $F=F_A$. It is not hard
to show that $\2D$ indeed is modular. 

If $Z_2(\2C)$ contains odd/Fermionic objects, i.e.\ objects with $\Theta_X=-\id_X$, the above
approach does not work. For a symmetric spherical category $\2C$ with non-trivial twists, one has a
generalization of Theorem \ref{theor-tann}, giving rise to a finite group $G$ together with an element 
$k\in Z(G)$ of order two. (Such a pair $(G,k)$ is occasionally called a super-group.)
One still has an equivalence $\2C\simeq\Rep\,G$ of fusion categories, but
the braidings of $\2C$ and $\Rep\,G$ are related by the Koszul-type rule
$c_{\2C}(X,Y)=\pm c_{\Rep\,G}(X,Y)$ for simple $X,Y$, where the minus sign applies when $X$ and $Y$
are both odd, cf.\ \cite[Section 7]{DR}. The regular representation of a super-group $(G,k)$  again is a
connected separable algebra $A\in\Rep(G,k)$, but it is only graded commutative. As a consequence,
when $\2C$ is a braided fusion category containing $\Rep(G,k)$ as a full subcategory, the module category 
${}_A\2C$ is not a $k$-linear tensor category but a tensor category enriched over the category
SVect$_k$ of super 
vector spaces. (Notice that, by contrast to tensor categories enriched over Vect, such a category is
{\it not} a tensor category, since the interchange law holds only up to signs: 
$(s\otimes t)\circ(s'\otimes t')=\pm(s\circ s')\otimes(t\circ t')$.)

However, if $(G,k)$ is a super-group, $\{e,k\}$ is a normal subgroup and if $H$ denotes the quotient 
group, $\Rep(G,k)$ contains the even category $\Rep\,H$ as a full fusion subcategory. If
$A\in\Rep\,H$ is the regular \'etale algebra of $H$, one finds ${}_A(\Rep(G,k))\cong$SVect, the
category of super-vector spaces, to wit the representation category of the super-group $(\{e,k\},k)$.
This shows that every pre-modular category $\2C$ admits a surjective braided functor
$F:\2C\rarr\2D$, where $\2D$ is `almost-modular':

\bdefin \label{def-almost} \cite{dno} A braided fusion category $\2C$ is called almost
non-degenerate if $Z_2(\2C)\simeq$SVect. (Equivalently, $Z_2(\2C)$ has precisely one non-trivial
simple object $X$ satisfying $X^2\cong\11$ and $\theta_X=-\id$.)
\edefin

Almost modular categories will briefly be mentioned again at the end of this review, when the recent
results of \cite{dno} will be touched upon.

The above considerations have an important generalization: Let $\2S\subset\2C$ be an arbitrary full
symmetric fusion subcategory of the pre-modular category $\2C$. For simplicity, we restrict to the
case where $\2S$ is even and thus equivalent (as a BTC) to the representation category of a finite
group $G$. Let $A\in\2S$ be the regular algebra. Since the latter contains all simple
objects of $\2S$ as direct summands, we have $A\in Z_2(\2C)$ if and only if $\2S\subset Z_2(\2C)$. 
If these equivalent conditions are satisfied, ${}_A\2C$ is braided and $F_A:\2C\rarr{}_A\2C$ is a
braided functor that trivializes the subcategory $\2S\subset\2C$. In that case, one finds
$Z_2({}_A\2C)\simeq {}_A(Z_2(\2C))$, which is trivial if and only if $\2S=Z_2(\2C)$, recovering the
previous result about modularization.

However, it is interesting to drop the requirement $\2S\subset Z_2(\2C)$. Independently of this
assumption, one finds that ${}_A\2C$ is a fusion category and $F_A$ a surjective tensor
functor. Furthermore, the group $G$ acts on the module category ${}_A\2C$ (by monoidal
self-equivalences)  and one has
$({}_A\2C)^G\simeq\2C$. When $\2S\not\subset Z_2(\2C)$, there exists no braiding on ${}_A\2C$ rendering
$F_A$ braided, but we have seen that there is a braided functor $\widehat{F_A}:\2C\rarr Z_1({}_A\2C)$.
In this specific situation, one can prove more: There is a $G$-grading on ${}_A\2C$, i.e.\ a map
$\del$ from the class of simple objects to $G$, constant on isomorphism classes and satisfying 
$\del(X\otimes Y)=\del X\cdot\del Y$. (As a consequence, $\del\11=e$ and 
$\del\ol{X}=(\del X)^{-1}$.) If $\2D$ is a $G$-graded category and $g\in G$, we denote by $\2D_g$
the full subcategory whose objects are direct sums of simple objects $X$ with $\del X=g$. Now we have
\[ ({}_A\2C)_e={}_A(C_\2C(\2S))={}_A\2C^0.\]
The action of $G$ on ${}_A\2C$ and the $G$-grading are connected by the identity
$\del(gX)=g(\del X)g^{-1}$, which is why ${}_A\2C$ is called a $G$-crossed category.
While ${}_A\2C$ does not admit a braiding (in the usual sense), it does admit a generalized braiding
that takes the grading and the $G$-action into account: For every $Y\in {}_A\2C$ and every 
$X\in({}_A\2C)_g$, there is an isomorphism $c_{X,Y}:X\otimes Y\rarr (gY)\otimes X$ satisfying
natural generalizations of the axioms for a braiding. Thus ${}_A\2C$ is a braided $G$-crossed
category. Conversely, if $\2D$ is a braided $G$-crossed fusion category then $\2C=\2D^G$ is an
ordinary braided fusion category containing $\Rep\,G$ as a full subcategory, and if $A$ is the
regular algebra in $\Rep\,G$ then ${}_A\2C\simeq\2D$. 

Most of these results are due to \cite{kir1,mue06,mue13} in the case of spherical fusion
categories. For a (somewhat) more extensive review than the one above, cf.\ \cite{mue-tur}. A much
longer discussion, including generalizations to not-necessarily-spherical fusion categories and
proofs pf precise 2-equivalences between categories of braided $G$-crossed fusion categories and
braided fusion categories containing $\Rep\,G$, cf.\ \cite{dgno2}.

Using the above results one can prove \cite{mue-tur} the following result concerning Conjecture
\ref{conj-embed}: 

\btheor \label{theor-4}
The following are equivalent:
\begin{itemize}
\item[(i)] Conjecture \ref{conj-embed} is true for every braided fusion category $\2C$ whose symmetric
center $Z_2(\2C)$ is even (and therefore equivalent to $\Rep\,G$ for a finite group $G$).
\item[(ii)] For every modular category $\2M$ acted upon by a finite group $G$ there is a braided
  crossed G-category $\2E$ with full $G$-spectrum and a $G$-equivariant equivalence $\2E_e\simeq\2M$.
\end{itemize}
\etheor

We close this section by pointing out that the above results have applications to the orbifold
construction in conformal field theory, cf.\ \cite[Section 6]{mue-tur} and references given there.


\section{The braided center of a fusion category}\label{s-z1}
As mentioned earlier, to every finite dimensional Hopf algebra $H$ one can associate \cite{drin} a
finite dimensional quasi-triangular Hopf algebra $(D(H),R)$, Drinfeld's `quantum double' of
$H$. (The construction is not restricted to finite dimensional algebras, but it becomes more
technical if $H$ is infinite dimensional and less relevant for the purposes of this review.) 
The case where $H$ is the group algebra of a finite group $G$ is denoted $D(G)$ and can be described
very explicitly. In particular, the braided category $D(G)$-$\Mod$ is modular \cite{ac}.

More generally, (i) semisimplicity of $D(H)$ is equivalent to (ii) semisimplicity of $H$ and of the
dual Hopf algebra $\widehat{H}$ and to (iii) $S_H^2=\id$ and $\dim H\ne 0$ in the ground field $k$.
Under these assumptions, $D(H)$-$\Mod$ is modular \cite{EG1}. (Cf.\ \cite[Appendix]{mue10} for an
alternative approach.)

Since the module category of a Hopf algebra $H$ satisfying the above conditions is a spherical
fusion category satisfying $\dim D(H)$-$\Mod=\dim_kH\ne 0$, the following results proven in
\cite{mue10} generalize those on $D(G)$ and $D(H)$:

\btheor \label{theor-double}
Let $k$ be an algebraically closed field and $\2C$ a $k$-linear semisimple spherical
category satisfying $\dim\2C\ne 0$. Then
\begin{itemize}
\item[(i)] $Z_1(\2C)$ is semisimple.
\item[(ii)] $Z_1(\2C)$ has a natural spherical structure inherited from $\2C$ and $\dim
  Z_1(\2C)=(\dim\2C)^2$. 
\item[(iii)] $Z_1(\2C)$ is non-degenerate, thus modular.
\item[(iv)] The Gauss sums (\ref{eq-gauss}) of $Z_1(\2C)$ are given by $\Omega^\pm(Z_1(\2C))=\dim\2C$.
\item[(v)] The forgetful functor $K: Z_1(\2C)\rarr\2C$ has a two-sided adjoint.
\item[(vi)] If $\2C$ already is modular, then the braided tensor functor 
$H:\2C\boxtimes\widetilde{\2C}\rarr Z_1(\2C)$ is an equivalence.
\item[(vii)] If $\2C_1,\2C_2$ satisfy the above assumptions and $\2C_1\approx\2C_2$ (monoidal
  Morita equivalence of \cite{mue09}) then $Z_1(\2C_1)\simeq Z_1(\2C_2)$ (braided equivalence).
\end{itemize}
\etheor

\brem 1. These results have been generalized to not necessarily spherical fusion categories,
cf.\ \cite{EO1}. Cf.\ also \cite{BV} for a more conceptual approach in terms of Hopf monads.

2. (vii) is an easy consequence of the definitions and a result in \cite{schauen}. The converse of
(vii) is also true, cf. \cite[Theorem 3.1]{eno2}. 

3. In view of (vi) and the considerations in Section \ref{ss-Z1}, one could take the equivalence
$Z_1(\2C)\simeq\2C\boxtimes\widetilde{\2C}$ as alternative definition of
modularity/non-degeneracy.

4. Statement (iii) means that $Z_2(Z_1(\2C))$ is trivial for any spherical fusion category
$\2C$. This should be compared with the other results of the type `the center of a center is 
trivial' mentioned in Remark \ref{rem-Z1}. It is tempting to conjecture that this holds more
generally in the context of centers in higher category theory.

5. Corollary (vi) implies  that every modular category $\2C$
arises as a direct factor of the braided center of some fusion category $\2D$: Just take
$\2D=\2C$. This is interesting since the braiding of $\2C$ is not used in defining
$Z_1(\2C)$. However, it seems pointless to reduce the classification of modular categories to the
classification of fusion categories, since there are many more of the latter and there is no hope of
classification. A more promising approach to `classifying' modular categories will be discussed in
the last section.

6. By (i)-(iii), the braided center construction gives rise to many modular categories. However, not
every modular category $\2C$ is equivalent to some $Z_1(\2D)$. This follows already from (iv) and the
fact that there are modular categories whose two Gauss sums are not equal. A criterion for
recognizing whether a modular category is of the form $Z_1(\2D)$ will be given below. Even when a modular
category $\2C$ does not satisfy this, one can often find fusion categories $\2D$ smaller than $\2C$ such that
$Z_1(\2D)$ has  $\2C$ as a direct factor.
\erem

We now turn to the question of recognizing the modular categories that are of the form $Z_1(\2D)$
for $\2D$ fusion, which has been solved quite recently, cf.\ \cite{dgno2,dmno}. As mentioned in
Section \ref{ss-algebras2}, a commutative (\'etale, connected) algebra $A$ in a braided fusion
category $\2C$ gives rise to a braided tensor functor $F_A:\2C\rarr Z_1({}_A\2C)$. Under very weak
conditions $F_A$ is faithful. In general, $F_A$ need not be full, but it is so when
$\2C$ is non-degenerate. This can be shown either by direct -- and tedious -- computation of
$\Hom_{Z_1(\2C)}(F_A(X),F_A(Y))$ or by invoking Proposition \ref{prop-faithful}.
Thus if $\2C$ is non-degenerate, $F_A:\2C\rarr Z_1({}_A\2C)$ is an embedding of braided fusion
categories, which by Theorem \ref{theor-factoriz}.(iv) gives rise to a direct factorization. The
complementary factor $C_{Z_1(\2C)}(F_A(\2C))'$ can be identified using the result of Schauenburg
mentioned in the last paragraph of Section \ref{ss-algebras2}.
Recall that if $\2C$ is non-degenerate braided fusion then we have the braided equivalence 
$\2C\boxtimes\widetilde{\2C}\stackrel{F_1\boxtimes F_2}{\longrightarrow} Z_1(\2C)$. 
If $A\in\2C$ is a commutative algebra then $B=F_2(A)\in Z_1(\2C)$ is a commutative algebra and
$\underline{B}=A$. The equivalence $Z_1(\2C)\rarr\2C\boxtimes\widetilde{\2C}$ maps $B$ to
$\11\boxtimes A$. Combining with Schauenburg's result, we have
\[ Z_1({}_A\2C)=Z_1({}_{\underline{B}}\2C)   \simeq {}_BZ_1(\2C)^0 
   \simeq{}_{(\11\boxtimes A)}(\2C\boxtimes\widetilde{\2C})^0 = \2C\boxtimes{}_A\widetilde{\2C}^0
= \2C\boxtimes\widetilde{{}_A\2C^0}.  \]
Thus, if $\2C$ is non-degenerate braided and $A\in\2C$ a connected \'etale algebra, there is a
braided equivalence
\be Z_1({}_A\2C)\simeq \2C\boxtimes\widetilde{{}_A\2C^0}. \label{eq-z1}\ee
Since one can prove, cf.\ \cite{dmno}, that $\2C=Z_1(\2D)$ contains a connected \'etale
algebra such that ${}_A\2C^0$ is trivial, one arrives at the following characterization of Drinfeld centers of fusion categories:

\btheor A non-degenerate braided fusion category $\2C$ is equivalent to $Z_1(\2D)$, where $\2D$ is a
fusion category, if and only if there is a connected \'etale algebra $A$ in $\2C$ such that
${}_A\2C^0$ is trivial. In this case, one can take $\2D={}_A\2C$.
\etheor

In Section \ref{s-modularization}, we defined fusion categories graded by a finite group. One can
ask how a $G$-grading on a fusion category $\2C$ is reflected in the center $Z_1(\2C)$. This was
clarified in \cite{GNN}, where the following is proven:

\btheor \cite[Theorem 3.5]{GNN} Let $\2C$ be $G$-graded with degree zero component $\2C_e$. Then the
relative Drinfeld 
center $Z_1(\2C,\2C_e)$ mentioned in Remark \ref{rem-Z1}.2 (which is monoidal but not braided) has a
natural structure of braided $G$-crossed category, and there is an equivalence
\[ Z_1(\2C,\2C_e)^G\simeq Z_1(\2C) \]
of braided categories.
\etheor

The interest of this theorem derives from the fact that the relative center $Z_1(\2C,\2C_e)^G$ may
be easier to determine than the full $Z_1(\2C)$.

We close this section with an important application of the braided center $Z_1$ and of Theorem
\ref{theor-double} to topology. Since $Z_1(\2C)$ is modular when $\2C$ is fusion, it gives
rise to a Reshetikhin-Turaev TQFT \cite{rt2,turaev1,turaev}. It is natural to ask whether there is a
more direct construction of this TQFT in terms of the spherical category $\2C$. In fact, shortly
after the Reshetikhin-Turaev construction, Turaev and Viro \cite{TV,turaev} proposed a construction
of $2+1$-dimensional TQFTs in terms of triangulations and `state-sums' rather than surgery. While
being fundamentally different from the RT-approach, the TV construction still required a modular category as
input. It was realized by various authors that the construction of a state-sum TQFT actually does
not require a braiding and that a spherical fusion category suffices as input datum, cf.\ in
particular \cite{BW2}. (The same observation was also made by Ocneanu and by S.~Gelfand and
Kazhdan.) This made it natural to conjecture that the state-sum TQFT of \cite{BW2} associated with a
spherical fusion category $\2C$ is isomorphic to the surgery TQFT associated with the modular
category $Z_1(\2C)$. This conjecture was proven in 2010, independently by Turaev and Virelizier
\cite{TurVir}, based on extensive previous work by Brugui\`eres and Virelizier \cite{BV,BLV}, and by
Balsam and Kirillov \cite{KB,balsam}.


\section{The Witt group of modular categories}
The results in this subsection are from \cite{dmno}. They are motivated by the desire to `classify'
modular categories (or non-degenerate braided fusion categories). This is a rather hopeless project
since, by Theorem \ref{theor-double}, $Z_1(\2D)$ is modular whenever $\2D$ is a spherical fusion
category and since there is little hope of classifying fusion categories. (Recall that, e.g., every
semisimple Hopf algebra gives rise to a fusion category.) The fact that Morita equivalent fusion
categories have equivalent Drinfeld centers reduces the problem only marginally.

This leads to the idea of considering categories of the form $Z_1(\2D)$ with $\2D$ fusion as
`trivial' and of classifying modular categories (or non-degenerate braided fusion
categories) `up to centers'. The following definition provides a rigorous way of doing this.

\bdefin Two non-degenerate braided fusion categories $\2C_1,\2C_2$ are called Witt equivalent if
there are fusion categories $\2D_1,\2D_2$ such that there is a braided equivalence
$\2C_1\boxtimes Z_1(\2D_1)\simeq \2C_2\boxtimes Z_1(\2D_2)$. 
\edefin

Witt equivalence obviously is coarser than braided equivalence, and it is not hard to show that it
is an equivalence relation. In fact, the Witt classes form a set $W_M$ that actually is countable.
Denoting the Witt equivalence class of $\2C$ by $[\2C]$, $W_M$ becomes an abelian monoid via
$[\2C_1]\cdot[\2C_2]:=[\2C_1\boxtimes\2C_2]$ with unit $\11_{W_M}=[\mbox{Vect}]$. Up to this point,
analogous results hold for the set of braided equivalence classes of non-degenerate braided fusion
categories, of which $W_M$ is the quotient monoid under the identification
$[Z_1(\2D)]=[\mbox{Vect}]$ for each fusion category $\2D$. But $W_M$ has one crucial additional
property: It is a group. Namely, defining $[\2C]^{-1}:=[\widetilde{\2C}]$, we have 
\[ [\2C]\cdot[\2C]^{-1}=[\2C]\cdot[\widetilde{\2C}]=[\2C\boxtimes\widetilde{\2C}]=[Z_1(\2C)]=\11_{W_M},\]
where the penultimate identity crucially depends on (vii) of Theorem \ref{theor-double}. Therefore,
$W_M$ is called the Witt group. (Actually, there are three Witt groups, defined in terms of
non-degenerate braided fusion categories, modular categories and unitary modular categories,
respectively.) 

By (\ref{eq-z1}), an \'etale algebra $A\in\2C$ gives rise to a braided equivalence 
$Z_1({}_A\2C)\simeq\2C\boxtimes\widetilde{{}_A\2C^0}$ and therefore to the identity
$[\2C]=[{}_A\2C^0]$ in the Witt group. Using the fact that $Z_1(\2D)$, where $\2D$ is fusion,
contains an \'etale algebra $A$ such that ${}_A\2C^0$ is trivial, the following is not hard to show,
cf.\ \cite{dmno}:

\btheor\label{theor-dmno}
Let $\2C_1, \2C_2$ be non-degenerate braided fusion categories. Then the following are equivalent:
\begin{itemize}
\item[(i)] $[\2C_1]=[\2C_2]$, i.e.\ $\2C_1$ and $\2C_2$ are Witt equivalent.
\item[(ii)] There is a fusion category $\2D$ such that 
$\2C_1\boxtimes\widetilde{\2C_2}\simeq Z_1(\2D)$. 
\item[(iii)] There is a connected \'etale algebra $A\in\2C_1\boxtimes\widetilde{\2C_2}$ such that 
${}_A(\2C_1\boxtimes\widetilde{\2C_2})^0$ is trivial.
\item[(iv)] There exist a non-degenerate braided fusion category $\2C$, connected \'etale algebras
  $A_1,A_2\in\2C$ and braided equivalences $\2C_1\simeq {}_{A_1}\2C^0,\ \2C_2\simeq {}_{A_2}\2C^0$. 
\item[(v)] There exist connected \'etale algebras $A_1\in\2C_1,A_2\in\2C_2$ and a braided equivalence
${}_{A_1}\2C_1^0\simeq{}_{A_2}\2C_2^0$.
\end{itemize}
\etheor

This shows that Witt equivalence could have been defined in terms of dyslectic module categories
instead of invoking the braided center $Z_1$. This latter approach has a `physical' interpretation:
Consider a rational chiral conformal field theory $\2A$, either as a ($C_2$-cofinite) vertex operator
algebra or in terms of von Neumann algebras indexed by intervals on $S^1$, as e.g.\ in
\cite{wass,xu,klm}. As mentioned earlier, in both settings there is a proof of modularity of the
representation category $\Rep\2A$. Furthermore, in both settings, there is a notion of `finite
extension' (or conformal extension)
and one can prove that the finite extensions $\2B\supset\2A$ are classified by the connected \'etale
algebras $A\in\Rep\2A$ in such way that $\Rep\2B={}_A(\Rep\2A)^0$ when the extension $\2B\supset\2A$ 
corresponds to the algebra $A\in\Rep\2A$. (Cf.\ \cite{ko} and 
\cite{mue-prag} for proof sketches.) This fact implies that we have $[\Rep\,\2B]=[\Rep\,\2A]$ for any
finite extension $\2B\supset\2A$ of rational chiral CFTs.

A (not very precise) folk conjecture in conformal field theory, cf.\ e.g.\ \cite{ms2}, states that
every modular category (to the extent that it is realized by a CFT) can be obtained from the modular
categories arising from WZW models combined and a certain set of `constructions' (like orbifold and
coset constructions). Now, the WZW categories coincide with the representation categories of quantum
groups at root-of-unity deformation parameter. Thus if one accepts that the above constructions
amount to passing to finite index subtheories and extensions, one arrives at the following
mathematical formulation of the Moore-Seiberg conjecture:

\bconj The Witt group $W_M$ is generated by the classes $[\2C(\6g,q)]$ of the quantum group
categories $\2C(\6g,q)$, where $\6g$ is a simple Lie algebra and $q$ a root of unity.
\econj

The only evidence for the conjecture so far is that there are no counterexamples! While there are
fusion categories that are `exotic' in the sense of having no (known) connection with finite group
theory or Lie theory, no modular categories are known that are `genuinely exotic' in the sense of
not being (related to) Drinfeld centers of exotic fusion categories. However, the existing
classification of conformal extensions provides a large and presumably complete set of relations in
the subgroup of the Witt group generated by the classes $[\2C(\6g,q)]$. While the full group $W_M$
is not understood, a close relative, to wit the Witt group of almost non-degenerate braided fusion
categories, has been computed recently, cf.\ \cite{dno}.

The circle of ideas around Witt equivalence is also relevant for the construction of two-dimensional
CFTs from a pair of chiral (`one-dimensional') CFTs. The relevant mathematical structure seems to be
the following:

\bdefin \label{def-inv}
Let $\2C_1,\2C_2$ be modular categories. A modular invariant for $(\2C_1,\2C_2)$ is a triple
$(A_1,A_2,E)$, where $A_1\in\2C_1,A_2\in\2C_2$ are connected \'etale algebras and
$E:{}_{A_1}\2C_1^0\rarr{}_{A_2}\2C_2^0$ is a braided equivalence.
\edefin

In view of Theorem \ref{theor-dmno}, it is clear that a modular invariant for $(\2C_1,\2C_2)$ exists
if and only if $\2C_1$ and $\2C_2$ are Witt equivalent. But more can be said (stated in
\cite{mue-prag} and proven in \cite{dno}):

\bprop If $\2C_1,\2C_2$ are modular, there is a bijection (modulo natural equivalence relations)
between modular invariants $(A_1,A_2,E)$ 
for $(\2C_1,\2C_2)$ and connected \'etale algebras $A\in\2C_1\boxtimes\widetilde{\2C_2}$ such that 
${}_A(\2C_1\boxtimes\widetilde{\2C_2})^0$ is trivial ($\Leftrightarrow\ d(A)=\sqrt{\dim\2C_1\cdot\dim\2C_2}$).
\eprop

A related result proven in \cite{ffrs} involves non-commutative algebras in a modular category. The
fact that an algebra over a field has a center, which is a commutative algebra, generalizes to
braided spherical categories. But since the definition of the center of an algebra $A$ in a braided
category $\2C$ involves the braiding, there will actually be two centers $Z_L(A), Z_R(A)$, depending
on the use of $c$ or $\widetilde{c}$. One finds

\btheor \cite{ffrs} Let $A$ be a separable connected algebra in a modular category $\2C$. Then the
centers 
$Z_L(A), Z_R(A)$ are connected \'etale and there is an equivalence 
\[ E: {}_{Z_L(A)}\2C^0 \rarr {}_{Z_R(A)}\2C^0 \]
of braided categories. (Thus $(Z_L(A),Z_R(A),E)$ is a modular invariant in the sense of Definition
\ref{def-inv}.) 

Every modular invariant arises in this way from a separable connected algebra $A$ in $\2C$, cf.\ \cite{KR}.
\etheor

\brem 1. The categorical constructions in \cite{ffrs} were inspired by analogous constructions in an
operator algebraic context, cf.\ \cite{BEK}, and conjectures in \cite{ostrik}.

2. In the series of papers \cite{furus}, a construction of ``topological two-dimensional CFTs'' was
given taking a modular category and a separable connected algebra in it as a starting point. (The
quotation marks refer to the fact that a CFT is more than a TQFT: It involves infinite dimensional
Hilbert spaces, trace class operators, analytic characters, etc.)
\erem


\end{document}